\documentclass[12pt]{article}

\usepackage[english]{babel}

\usepackage[pdftex,paper=a4paper,portrait=true,textwidth=450pt,textheight=675pt,tmargin=3cm,marginratio=1:1]{geometry}

\usepackage{amsfonts}

\usepackage[dvips]{graphics}

\usepackage{amsmath}

\usepackage{amsthm}

\usepackage{amssymb}

\usepackage{eufrak}

\usepackage{cancel}

\renewcommand\L{\mathcal L}

\renewcommand\O{\mathcal O}
\newcommand\PP{\mathbb P}

\newcommand\D{\mathcal D}
\newcommand\bsdelta{\boldsymbol \delta}
\newcommand\bschi{\boldsymbol \chi}
\newcommand\bslambda{\boldsymbol \lambda}

\newcommand\F{\mathcal F}
\newcommand\G{\mathcal G}

\newcommand\M{\mathcal M}
\newcommand\N{\mathcal N}

\newcommand\C{\mathbb C}

\newcommand\FF{\mathbb F}

\newcommand\Q{\mathbb Q}
\newcommand\cQ{\mathcal Q}

\newcommand\ccR{\mathcal R}

\newcommand\cT{\mathcal T}

\newcommand\X{\mathcal X}

\newcommand\Z{\mathbb Z}

\newcommand\Into{\ar@{^(->}[r]<-.3ex>}

\newcommand\rk{\operatorname{rk}}

\newcommand\ch{\operatorname{ch}}

\newcommand\Pic{\operatorname{Pic}}

\newcommand\Hilb{\operatorname{Hilb}}

\newcommand\beq[1]{\begin{equation}\label{#1}}
\newcommand\eeq{\end{equation}}
\newcommand\beqa{\begin{eqnarray*}}
\newcommand\eeqa{\end{eqnarray*}}

\input xy

\xyoption{all}

\newtheorem{theorem}{Theorem}[section]

\newtheorem{corollary}[theorem]{Corollary}

\newtheorem{proposition}[theorem]{Proposition}

\theoremstyle{definition}

\newtheorem{definition}[theorem]{Definition}

\newtheorem{remark}[theorem]{Remark}

\newtheorem*{acknowledgements}{Acknowledgements}

\theoremstyle{property}

\makeatletter
\newcommand{\contraction}[5][1ex]{%
  \mathchoice
    {\contraction@\displaystyle{#2}{#3}{#4}{#5}{#1}}%
    {\contraction@\textstyle{#2}{#3}{#4}{#5}{#1}}%
    {\contraction@\scriptstyle{#2}{#3}{#4}{#5}{#1}}%
    {\contraction@\scriptscriptstyle{#2}{#3}{#4}{#5}{#1}}}%
\newcommand{\contraction@}[6]{%
  \setbox0=\hbox{$#1#2$}%
  \setbox2=\hbox{$#1#3$}%
  \setbox4=\hbox{$#1#4$}%
  \setbox6=\hbox{$#1#5$}%
  \dimen0=\wd2%
  \advance\dimen0 by \wd6%
  \divide\dimen0 by 2%
  \advance\dimen0 by \wd4%
  \vbox{%
    \hbox to 0pt{%
      \kern \wd0%
      \kern 0.5\wd2%
      \contraction@@{\dimen0}{#6}%
      \hss}%
    \vskip 0.5ex
    \vskip\ht2}}

\newcommand{\contraction@@}[3][0.05em]{%
  \hbox{%
    \vrule width #1 height 0pt depth #3%
    \vrule width #2 height 0pt depth #1%
    \vrule width #1 height 0pt depth #3%
    \relax}}
\makeatother

\begin{document}
\title{Euler characteristics of moduli spaces of torsion free sheaves on toric surfaces}
\author{Martijn Kool\thanks{Pacific Institute for the Mathematical Sciences, University of British Columbia, 4176--2207 Main Mall, Vancouver, British Columbia, Canada, V6T 1Z4, {\tt{mkool@math.ubc.ca}}.}}
\maketitle
\begin{abstract}
Given a smooth toric variety $X$, the action of the torus $T$ lifts to the moduli space $\M$ of stable sheaves on $X$. Using the pioneering work of Klyacho, a fairly explicit combinatorial description of the fixed point locus $\M^T$ can be given (as shown by earlier work of the author). In this paper, we apply this description to the case of torsion free sheaves on a smooth toric surface $S$.  A general expression for the generating function of the Euler characteristics of such moduli spaces is obtained. The generating function is expressed in terms of Euler characteristics of certain moduli spaces of stable configurations of linear subspaces appearing in classical GIT. The expression holds for \emph{any} choice of $S$, polarization, rank, and first Chern class. Specializing to various examples allows us to compute some new as well as known generating functions. 
\end{abstract}

\section{Introduction}

The moduli space $\M$ of Gieseker stable\footnote{For the definition of Gieseker stability, see \cite[Def.~1.2.4]{HL}.} sheaves is a complicated object. It satisfies Murphy's Law, meaning every singularity type of finite type over $\mathbb{Z}$ appears on one of its components \cite{Vak}. Many geometrically interesting invariants are defined on components of this moduli space and their computation requires us to have some understanding of these components. Examples of invariants are motivic invariants such as Euler characteristic or (virtual) Poincar\'e polynomial. Another example is the Donaldson-Thomas invariants of a Calabi-Yau 3-fold. 

Let $X$ be a polarized\footnote{Recall that the notion of stability depends on the choice of polarization.} smooth projective toric variety\footnote{In this paper, we work with varieties, schemes, and stacks over ground field $\mathbb{C}$.} with torus $T$. The action of $T$ on $X$ lifts to $\M$. One can hope that this action facilitates explicit computation of invariants of $\M$ by reduction to the fixed point locus $\M^T \subset \M$. Based on ideas of Klyachko \cite{Kly1, Kly2, Kly3, Kly4}, the author gives a fairly explicit description of the fixed point locus $\M^T$ in \cite{Koo1}. In the case of $\mu$-stability\footnote{For the definition of $\mu$-stability, also known as slope or Mumford-Takemoto stability, see \cite[Def.~1.2.12]{HL}.} and reflexive sheaves, this description simplifies significantly \cite{Koo1}. In the present paper, we systematically specialize these ideas to the case $X=S$ is a toric surface. For applications to pure dimension 1 sheaves on toric surfaces, see \cite{Cho}, \cite{CM}, and \cite[Sect.~2.4]{Koo2}. 

Let $S$ be a smooth complete toric surface with polarization $H$. Denote by $\M_{S}^{H}(r,c_1,c_2)$ the moduli space of $\mu$-stable torsion free sheaves on $S$ with rank $r$ and Chern classes $c_1, c_2$. The main result of this paper is an expression for the generating function
\begin{equation} \label{genfun}
\sum_{c_2} e(\M_{S}^{H}(r,c_1,c_2)) q^{c_2},
\end{equation}
for \emph{any} $S$, $H$, $r$, $c_1$. Here $e(\cdot)$ denotes topological Euler characteristic. The expression is in terms of Euler characteristics of moduli spaces of stable configurations of linear subspaces\footnote{Configurations of linear subspaces and their moduli spaces are a classical topic in GIT. See \cite[Ch.~11]{Dol} for a discussion.} (Theorem \ref{ch. 2, sect. 3, thm. 1}). The expression can be further simplified in examples. The dependence on $H$ allows us to study wall-crossing phenomena in examples. Note that we compute Euler characteristics of moduli spaces of $\mu$-stable torsion free sheaves only, even when strictly $\mu$-semistable torsion free sheaves are present.  

This paper is organized as follows. In Section 2, we recall the relevant results from \cite{Koo1}. In Section 3, we give an explicit formula for the Chern character of an arbitrary $T$-equivariant locally free sheaf on $S$. Each torsion free sheaf on $S$ embeds in its double-dual, which is reflexive and hence locally free (because $\dim(S)=2$). Using the double-dual map, the generating function (\ref{genfun}) can be written as a product of an explicit 0-dimensional part times 
\begin{equation} \label{genfunrfl}
\sum_{c_2} e(\N_{S}^{H}(r,c_1,c_2)) q^{c_2},
\end{equation}
where $\N_{S}^{H}(r,c_1,c_2)$ is the moduli space of $\mu$-stable locally free sheaves on $S$ with rank $r$ and Chern classes $c_1, c_2$. This product structure was first pointed out by G\"ottsche and Yoshioka \cite[Prop.~3.1]{Got3}. The generating function (\ref{genfunrfl}) can be expressed explicitly in terms of Euler characteristics of moduli spaces of stable configurations of linear subspaces (Theorem \ref{ch. 2, sect. 3, thm. 1}). In Section 4, we apply the formula to various examples and compare to results in the literature. For rank 1, this gives the formula of Ellingsrud and Str{\o}mme \cite{ES} and G\"ottsche \cite{Got1}. Note that G\"ottsche's formula holds on any smooth complete surface. For rank 2 and $S = \PP^2$, we obtain a simple formula which we compare to work of Klyachko \cite{Kly4} and Yoshioka \cite{Yos}. For rank 2 and $S = \mathbb{P}^{1} \times \mathbb{P}^{1}$ or any Hirzebruch surface $\mathbb{F}_{a}$, we make the dependence on choice of ample divisor $H$ explicit. This allows us to study wall-crossing phenomena and compare to work of G\"ottsche \cite{Got2} and Joyce \cite{Joy2}. We perform various consistency checks. Finally, we compute\footnote{This example was considered independently around the same time by Weist using techniques of toric geometry and quivers \cite{Wei}.} an explicit expression for rank 3 and $S = \mathbb{P}^{2}$. We would like to point out that \cite{ES, Kly4} use torus localization, whereas \cite{Got1, Got2, Yos} use very different techniques namely the Weil Conjectures. Also \cite{Joy2} uses very different techniques namely his theory of wall-crossing for motivic invariants counting (semi)stable objects in an abelian category. 

Finally, we would like to point out some important related literature. In \cite{BPT}, Bruzzo, Poghossian, and Tanzini study moduli spaces of framed torsion free sheaves on Hirzebruch surfaces using localization techniques. Furthermore, after the appearance of the preprint version of this paper, Manschot addressed modularity of the rank 3 generating function on $S = \PP^2$. Using a blow-up formula to get from $\PP^2$ to $\FF_1$ and a wall-crossing computation on $\FF_1$, he computes an expression for the generating function in terms modular forms and indefinite theta functions \cite{Man1}. Further recent computations on rational and ruled surfaces can be found in \cite{Man2, Man3, Moz}. \\

\noindent \textbf{Notation.} Two pieces of notation. (1) We denote by $\mathrm{Gr}(k,n)$ the Grassmannian of $k$-dimensional subspaces $V \subset \C^{\oplus n}$. (2) Let $a, b \in \mathbb{Z}$ with $a \neq 0$. We write $a \ | \ b$ whenever $b = a k$ for some $k \in \mathbb{Z}$.

\noindent \begin{acknowledgements} Many of the guiding ideas of this paper come from Klyachko's wonderful preprint \cite{Kly4}. I would like to thank Tom Bridgeland, Frances Kirwan, Sven Meinhardt, Yinan Song, Bal\'azs Szendr\H{o}i, Yukinobu Toda and Richard Thomas for useful discussions and my supervisor Dominic Joyce for his continuous support. The author would also like to thank the referee very useful suggestions on improving the exposition of this paper. This paper is part of the author's D.Phil.~project funded by an EPSRC Studentship, which is part of EPSRC Grant EP/D077990/1. 
\end{acknowledgements}

\section{Moduli spaces of sheaves on toric varieties}

This section is a brief exposition of the main results of \cite{Kly2, Kly4, Per, Koo1}. We review Klyacho's and Perling's descriptions of $T$-equivariant coherent, torsion free, and reflexive sheaves on toric varieties. We also discuss Klyachko's formula for the Chern character of a $T$-equivariant torsion free sheaf.

\subsection{Equivariant sheaves on toric varieties}

Let $X$ be a smooth toric variety of dimension $d$ with torus $T$. Let $M = X(T)$ be the character group of $T$ (written additively) and denote its dual by $N$. Denote the natural pairing by $\langle \cdot, \cdot \rangle  : M \times N \rightarrow \Z$. Then $N$ is a rank $d$ lattice containing a fan\footnote{We always assume $\Delta$ contains cones of dimension $d$.} $\Delta$ and the data $(N,\Delta)$ completely describes $X$. We refer to Fulton's book \cite{Ful} for the general theory. We recall that there is a bijection between the cones $\sigma \in \Delta$ and the $T$-invariant affine open subsets $U_\sigma \subset X$. \\

\noindent \emph{The affine case}. Suppose $X = U_\sigma$. Let $S_{\sigma} = \{ m \in M \ : \ \langle m , \sigma \rangle \geq 0 \}$. This semi-group gives rise to an algebra $\C[S_\sigma]$, which is exactly the coordinate ring of $U_\sigma$. Therefore, quasi-coherent sheaves on $U_\sigma$ are the same as $\C[S_\sigma]$-modules. More precisely, the global section function gives an equivalence of categories
\[
H^0(\cdot) : \mathrm{Qco}(U_\sigma) \rightarrow \C[S_\sigma]\textrm{-}\mathrm{Mod}.
\]
Under this equivalence, coherent sheaves correspond to the finitely generated modules. It will not come as a surprise that this equivalence can be extended to an equivalence between the categories of $T$-equivariant quasi-coherent sheaves and $\C[S_\sigma]$-modules with regular $T$-action. For a $T$-equivariant quasi-coherent sheaf $(\F,\Phi)$ on $U_\sigma$, use the $T$-equivariant structure $\Phi$ to define a regular $T$-action on $H^0(\F)$. Since $T$ is diagonalizable, a $T$-action on $H^0(\F)$ is equivalent to a decomposition of $H^0(\F)$ into weight spaces
\[
H^0(\F) = \bigoplus_{m \in M} H^0(\F)_m.
\] 
Therefore $T$-equivariant quasi-coherent sheaves on $U_\sigma$ are nothing but $M$-graded $\C[S_\sigma]$-modules, i.e.~there exists an equivalence of categories
\[
H^0(\cdot) : \mathrm{Qco}^T(U_\sigma) \rightarrow \C[S_\sigma]\textrm{-}\mathrm{Mod}^{M\textrm{-graded}}.
\]
See \cite{Kan}, \cite{Per} for details. \\

\noindent \emph{Repackaging in terms of $\sigma$-families}. Following Perling \cite{Per}, we write the data of an $M$-graded $\C[S_\sigma]$-module in a slightly more explicit way. 
\begin{definition}[Perling]
For each $m,m' \in M$ we write $m \leq_\sigma m'$ when $m' - m \in S_\sigma$. A $\sigma$-family $\hat{F}$ consists of the following data: a collection of complex vector space $\{F_m\}_{m \in M}$ and linear maps $\{\chi_{m,m'} : F_m \rightarrow F_{m'}\}_{m \leq_\sigma m'}$ such that:
\begin{enumerate}
\item[(i)] $\chi_{m,m} = \mathrm{id}_{F_m}$,
\item[(ii)] $\chi_{m',m^{\prime \prime}} \circ \chi_{m,m'} = \chi_{m, m^{\prime \prime}}$ for all $m \leq_\sigma m' \leq_\sigma m''$.
\end{enumerate}
A morphism between $\sigma$-families $\hat{F}, \hat{G}$ is a collection $\hat{\phi}$ of linear maps $\{\phi_m : F_m \rightarrow G_m\}_{m \in M}$ commuting with the $\chi$'s. \hfill $\oslash$
\end{definition}
An $M$-graded module $F =\bigoplus_{m \in M} F_m$ gives rise to a $\sigma$-family as follows. We simply take $\{F_m\}_{m \in M}$ to be the collection of weight spaces. For each $m \leq_\sigma m'$ we have $m' - m \in S_{\sigma} \subset M$, so multiplication by the character $m' - m$ gives a linear map $F_m \rightarrow F_{m'}$. This gives an equivalence of categories \cite[Prop.~5.5]{Per}
\[
\C[S_\sigma]\textrm{-}\mathrm{Mod}^{M\textrm{-graded}} \rightarrow \sigma\textrm{-Families}.
\]
When $\sigma$ is a cone of maximal dimension $d$, we can choose an order of its rays $(\rho_1, \ldots, \rho_d)$ and choose a primitive generator $n_i$ of each ray $\rho_i$. By smoothness of $U_\sigma$, this gives a basis $(n_1, \ldots, n_d)$ of the lattice $N$. Denote the dual basis by $(m_1, \ldots, m_d)$. This choice induces an isomorphism $U_\sigma \cong \C^d$. Let $\hat{F}$ be a $\sigma$-family. Writing each $m \in M$ as $m = \sum_i \lambda_i m_i$, we define
\[
F(\lambda_1, \ldots, \lambda_d) := F_m.
\]
Moreover, multiplication by $\chi_{m, m+m_i}$ gives linear maps
\[
\chi_i(\lambda_1, \ldots, \lambda_d):= \chi_{m, m+m_i} : F(\lambda_1, \ldots, \lambda_{d}) \rightarrow  F(\lambda_1, \ldots, \lambda_{i-1},\lambda_i+1,\lambda_{i+1}, \ldots ,\lambda_d)
\]
satisfying the usual commutativity requirements. We note some important properties.
\begin{enumerate}
\item[(i)] Let $\F$ be a $T$-equivariant quasi-coherent sheaf with $\sigma$-family $\hat{F}$. Then $\F$ is coherent if  only if $\hat{F}$ has finitely many homogeneous generators. We call such $\sigma$-families finite \cite[Def.~5.10]{Per}. 
\item[(ii)] Let $\F$ be a $T$-equivariant coherent sheaf with $\sigma$-family $\hat{F}$. Then $\F$ is torsion free if  only if all maps $\{\chi_{m,m'}\}_{m \leq_\sigma m'}$ are injective. This can be seen by noting that a non-trivial kernel of some $\chi_{m,m'}$ would give rise to a lower dimensional $T$-equivariant subsheaf of $\F$, which violates torsion freeness (e.g.~see \cite[Prop.~2.8]{Koo1}).
\end{enumerate}

\noindent \emph{Equivariant torsion free sheaves}. Let $\F$ be an $T$-equivariant coherent sheaf on $X$. Let $\{\sigma_1, \ldots, \sigma_e\}$ be the cones of maximal dimension. Note that $e = e(X)$ is the number of $T$-fixed points of $X$, which is equal to the Euler characteristic of $X$. The open subsets $U_{\sigma_i} \cong \C^d$ provide a $T$-invariant open affine cover of $X$ and the restrictions $\F|_{U_{\sigma_i}}$ give us a collection of finite $\sigma$-families $\{\hat{F}^{\sigma_i}\}_{i =1, \ldots, e}$. Now suppose we are given \emph{any} collection of finite $\sigma$-families $\{\hat{F}^{\sigma_i}\}_{i =1, \ldots, e}$. When do they ``glue'' to an $T$-equivariant coherent sheaf on $X$? In this paper, we are only interested in the torsion free case, so we describe the answer in this case only\footnote{For gluing conditions for general $T$-equivariant coherent sheaves see \cite[Sect.~5.2]{Per}.}. As mentioned above, in the torsion free case all the maps $\chi_{m,m'}^{\sigma_i}$ between the weight spaces are injective. We can assume all these maps are actually inclusions\footnote{The precise statement is this. The category of $T$-equivariant torsion free sheaves on $U_{\sigma_i}$ is equivalent to the category of finite $\sigma_i$-families with all maps $\chi_{m,m'}^{\sigma_i}$ injective. This category is equivalent to its full subcategory of finite $\sigma_i$-families with all maps $\chi_{m,m'}^{\sigma_i}$ inclusions.}. 

We now describe the gluing conditions. For each $i=1, \ldots, e$, let $(\rho^{(i)}_1, \ldots, \rho^{(i)}_d)$ be an ordering of rays of $\sigma_i$. Fix any two $i,j$, then the intersection $\sigma_i \cap \sigma_j$ is a cone of some dimension $p$. Assume w.l.o.g.~that $\sigma_i \cap \sigma_j$ is spanned by the first $p$ rays among $(\rho^{(i)}_1, \ldots, \rho^{(i)}_d)$ and $(\rho^{(j)}_1, \ldots, \rho^{(j)}_d)$. Then the corresponding gluing condition is
\begin{equation} \label{glue}
F^{\sigma_i}(\lambda_1, \ldots, \lambda_p, \infty, \ldots, \infty) = F^{\sigma_j}(\lambda_1, \ldots, \lambda_p, \infty, \ldots, \infty), \forall \ \lambda_1, \ldots, \lambda_p \in \Z.
\end{equation}
This needs some explanation. For fixed $\lambda_1, \ldots, \lambda_p \in \Z$ consider 
\[
\{F^{\sigma_i}(\lambda_1, \ldots, \lambda_p, \mu_{p+1}, \ldots, \mu_d)\}_{\mu_{p+1}, \ldots, \mu_{d} \in \Z}.
\]
Since the $\sigma$-family $\hat{F}^{\sigma_i}$ is finite, these vector spaces stabilize for sufficiently large $\mu_{p+1}$, $\ldots$, $\mu_{d}$ and we denote the limit by $F^{\sigma_i}(\lambda_1, \ldots, \lambda_p, \infty, \ldots, \infty)$. Moreover, the vector spaces $F^{\sigma_i}(\lambda_1, \ldots, \lambda_d)$ form a multi-filtration of some limiting finite dimensional vector space $F^{\sigma_i}(\infty, \ldots, \infty)$ of dimension $\rk(\F)$. The idea is that the left hand side of (\ref{glue}) is the $\sigma$-family of $\F|_{U_{\sigma_i}}$ restricted to $U_{\sigma_i} \cap U_{\sigma_{j}}$ and the right hand side is the $\sigma$-family of $\F|_{U_{\sigma_j}}$ restricted to $U_{\sigma_i} \cap U_{\sigma_{j}}$. This description of $T$-equivariant torsion free sheaves is originally due to Klyachko \cite{Kly2, Kly4}. We summarize:
\begin{theorem}[Klyachko] \label{Kly}
Let $X$ be a smooth toric variety described by a fan $\Delta$ in a lattice $N$ of rank $d$. Let $\{\sigma_1, \ldots, \sigma_e\}$ be the cones of maximal dimension. For each $i=1, \ldots, e$, let $(\rho^{(i)}_1, \ldots, \rho^{(i)}_d)$ be an ordering of the rays of $\sigma_i$. The category of $T$-equivariant torsion free sheaves on $X$ is equivalent to a category $\cT$ which can be described as follows. The objects of $\cT$ are collections of finite $\sigma$-families $\{\hat{F}^{\sigma_i} \}_{i = 1, \ldots, e}$, with all maps $\chi_{m,m'}^{\sigma_i}$ inclusions, satisfying the following gluing condition. For any two $i,j$, $\sigma_i \cap \sigma_j$ is a cone of some dimension $p$. Assume w.l.o.g.~that $\sigma_i \cap \sigma_j$ is spanned by the first $p$ rays among both $(\rho^{(i)}_1, \ldots, \rho^{(i)}_d)$ and $(\rho^{(j)}_1, \ldots, \rho^{(j)}_d)$. Then $\hat{F}^{\sigma_i}$, $\hat{F}^{\sigma_j}$ satisfy\footnote{It should be clear how the gluing conditions read when the rays of $\sigma_i \cap \sigma_j$ do not necessarily correspond to the first $p$ rays of $\sigma_i$ and $\sigma_j$.} (\ref{glue}). The maps of $\cT$ are collections of maps of $\sigma$-families $\{\hat{\phi}^{\sigma_i} : \hat{F}^{\sigma_i} \rightarrow \hat{G}^{\sigma_i} \}_{i = 1, \ldots, e}$ such that for each $i,j$ as above$^{10}$ 
\[
\phi^{\sigma_i}(\lambda_1, \ldots, \lambda_p, \infty, \ldots, \infty) = \phi^{\sigma_j}(\lambda_1, \ldots, \lambda_p, \infty, \ldots, \infty), \ \forall \ \lambda_1, \ldots, \lambda_p \in \Z.
\]
\end{theorem}

Although the description in this theorem is not entirely coordinate invariant, the only choice we made is an ordering of the rays of each cone $\sigma_i$ of maximal dimension. For an extension of this theorem to any $T$-equivariant pure sheaves, see \cite[Sect.~2]{Koo1}. \\

\noindent \emph{Equivariant reflexive sheaves}. Let $(\cdot)^* = \mathcal{H}{\it{om}}(\cdot, \O_X)$. A coherent sheaf $\F$ on $X$ is called reflexive if the natural morphism $\F \rightarrow \F^{**}$ is an isomorphism. A $T$-equivariant reflexive sheaf on $X$ is $T$-equivariant torsion free. However, $T$-equivariant reflexive sheaves have a simpler description than $T$-equivariant torsion free sheaves. The reason is that reflexive sheaves are fully determined by their behaviour off any codimension $\geq 2$ closed subset \cite[Prop.~1.6]{Har}. In particular, a reflexive sheaf on a $T$-invariant affine open subset $U_{\sigma_i} \cong \C^d$ is fully determined by its restriction to the complement of the union of all codimension 2 coordinate hyperplanes
\[
(\C \times \C^* \times \cdots \times \C^*) \cup (\C^* \times \C \times \C^* \times \cdots \times \C^*) \cup \cdots \cup (\C^* \times \cdots \times \C^* \times \C).
\]
The restrictions to the components of this union are easy to describe. We give the final result:

Let $\Delta(1)$ be the collection of rays of the fan $\Delta$ of $X$. We introduce a category $\ccR$. Its objects are collections of vector spaces $\{V^\rho(\lambda)\}_{\rho \in \Delta(1), \lambda \in \Z}$ which form flags
\[
\cdots \subset V^\rho(\lambda-1) \subset V^\rho(\lambda) \subset V^\rho(\lambda+1) \subset \cdots.
\] 
We require these flags to be finite meaning $V^\rho(\lambda) = 0$ for $\lambda \ll 0$. They are also required to be full meaning $V^\rho(\lambda) = V^\rho(\lambda+1)$ for $\lambda \gg 0$. We denote the limiting vector space by $V^\rho(\infty)$. The maps in the category $\ccR$ are the obvious: linear maps between the limiting vector spaces preserving the flags. There is a natural fully faithful functor $\ccR \rightarrow \cT$ defined as follows. As before, denote the cones of $\Delta$ of maximal dimension by $\sigma_1, \ldots, \sigma_e$. For each $i=1,\ldots, e$, let $(\rho^{(i)}_1, \ldots, \rho^{(i)}_d)$ be an ordering of rays of $\sigma_i$. Then we map  $\{V^\rho(\lambda)\}_{\rho \in \Delta(1), \lambda \in \Z}$ to the following collection of finite $\sigma$-families
\begin{equation*} 
F^{\sigma_i}(\lambda_1, \ldots, \lambda_d) := V^{\rho^{(i)}_{1}}(\lambda_1) \cap \cdots \cap V^{\rho^{(i)}_{d}}(\lambda_d), \ \forall \lambda_1, \ldots, \lambda_d \in \Z.
\end{equation*}
Under the equivalence of categories of Theorem \ref{Kly}, the $T$-equivariant reflexive sheaves on $X$ correspond to the elements of the image of $\ccR \rightarrow \cT$ \cite{Kly1, Kly2}, \cite[Thm.~5.19]{Per}. From the fact that rank 1 reflexive sheaves are line bundles, one easily deduces that the $T$-equivariant Picard group $\Pic^T(X)$ is isomorphic to $\Z^{\#\Delta(1)}$.

\subsection{Moduli spaces of equivariant sheaves}

Theorem \ref{Kly} allows one to construct explicit moduli spaces of $T$-equivariant torsion free sheaves. A natural topological invariant of a $T$-equivariant sheaf is its characteristic function \cite[Def.~3.1]{Koo1}. Again, in this section we only consider the torsion free case\footnote{Large parts of this section hold for $T$-equivariant pure sheaves in  general \cite{Koo1}.}.
\begin{definition} \label{ch. 1, sect. 3, def. 1}
Let the notation be as in Theorem \ref{Kly}. Let $\F$ be a $T$-equivariant torsion free sheaf on $X$, then the \emph{characteristic function} $\bschi_{\F}$ of $\F$ is 
\begin{align} 
\begin{split}
&\bschi_{\F} : M \longrightarrow \mathbb{Z}^{e}, \\
&\bschi_{\F}(m) = (\chi_{\F}^{\sigma_{1}}(m), \ldots, \chi_{\F}^{\sigma_{e}}(m)) = (\mathrm{dim}(F_{m}^{\sigma_{1}}), \ldots, \mathrm{dim}(F_{m}^{\sigma_{e}})). \nonumber
\end{split}
\end{align}
We denote the set of all characteristic functions by $\mathcal{X}$. \hfill $\oslash$
\end{definition}

Given a $T$-equivariant $S$-flat family of coherent sheaves, it is not hard to see that characteristic functions are locally constant on the base $S$ \cite[Prop.~3.2]{Koo1}. This makes it a good topological invariant. Moreover, it is finer than Hilbert polynomial. More precisely, fixing a polarization on $X$, any two $T$-equivariant torsion free sheaves on $X$ with the same characteristic function $\bschi$ have the same Hilbert polynomial \cite[Prop. 3.14]{Koo1}. We refer to this polynomial as the Hilbert polynomial determined by $\bschi$. For a fixed Hilbert polynomial $P$, we denote by $\mathcal{X}_{P} \subset \X$ the set of characteristic functions which determine the Hilbert polynomial $P$. 

For any $\bschi \in \X$, one can now define moduli functors 
\begin{align*} 
\underline{\mathcal{M}}_{\bschi}^{ss} &: (Sch/\mathbb{C})^{o} \longrightarrow Sets \\
\underline{\mathcal{M}}_{\bschi}^{s} &: (Sch/\mathbb{C})^{o} \longrightarrow Sets 
\end{align*}   
of $T$-equivariant flat families\footnote{As usual, two such families $\mathcal{F}_{1}$, $\mathcal{F}_{2}$ are identified if there exists a line bundle $L$ on $S$ and a $T$-equivariant isomorphism $\mathcal{F}_{1} \cong \mathcal{F}_{2} \otimes p_{S}^{*}L$. See \cite[Sect.~3.1]{Koo1} for details.} with fibres Gieseker semistable (respectively geometrically Gieseker stable) $T$-equivariant torsion free sheaves on $X$ with characteristic function $\bschi$.

Using Theorem \ref{Kly}, it is a straight-forward exercise in GIT to define candidate schemes $\mathcal{M}_{\bschi}^{\tau, ss}$, $\mathcal{M}_{\bschi}^{\tau, s}$ corepresenting these functors. One takes certain closed subschemes of products of Grassmannians (describing the multi-filtrations of Theorem \ref{Kly}) and considers the natural $G = \mathrm{SL}(r,\C)$ action on it. Here $r = \chi^{\sigma_{1}}(\infty, \ldots, \infty) = \cdots = \chi^{\sigma_{e}}(\infty, \ldots, \infty)$ is the dimension of the limiting vector space. Then two objects are $T$-equivariantly isomorphic if and only if the corresponding points lie in the same $G$-orbit. The hard part is to find a $G$-equivariant line bundle which reproduces Gieseker stability. Such $G$-equivariant line bundles are constructed in \cite[Thm.~3.21]{Koo1}.
\begin{theorem} [{{\cite[Thm.~3.12]{Koo1}}}] \label{ch. 1, sect. 3, thm. 2} 
Let $X$ be a polarized smooth projective toric variety and let $\bschi \in \mathcal{X}$. Then $\underline{\mathcal{M}}_{\bschi}^{ss}$ is corepresented by a projective scheme $\mathcal{M}^{ss}_{\bschi}$ explicitly constructed using GIT in \cite[Sect.~3.3]{Koo1}. Moreover, there is an open subset $\mathcal{M}_{\bschi}^{s} \subset \mathcal{M}_{\bschi}^{ss}$ such that $\underline{\mathcal{M}}_{\bschi}^{s}$ is corepresented by $\mathcal{M}^{s}_{\bschi}$ and $\mathcal{M}^{s}_{\bschi}$ is a coarse moduli space. 
\end{theorem}

The construction of the moduli spaces $\mathcal{M}^{ss}_{\bschi}$, $\mathcal{M}^{s}_{\bschi}$ simplifies considerably if one replaces ``torsion free'' by ``reflexive'' and ``Gieseker stable'' by ``$\mu$-stable'' \cite[Sect.~4.4]{Koo1}. Denote by $\mathcal{X}^{\mathrm{refl}} \subset \mathcal{X}$ be the subset of characteristic functions of $T$-equivariant reflexive sheaves on $X$. For any $\bschi \in \mathcal{X}^{\mathrm{refl}}$, define moduli functors
\begin{align*} 
\underline{\mathcal{N}}_{\bschi}^{\mu ss} &: (Sch/\mathbb{C})^{o} \longrightarrow Sets \\
\underline{\mathcal{N}}_{\bschi}^{\mu s} &: (Sch/\mathbb{C})^{o} \longrightarrow Sets 
\end{align*}  
of $T$-equivariant $S$-flat families$^{12}$ with fibres $\mu$-semistable (resp.~geometrically $\mu$-stable) $T$-equivariant reflexive sheaves on $X$ with characteristic function $\bschi$. Again, straightforward use of GIT yields candidate schemes $\mathcal{N}_{\bschi}^{\mu ss}$, $\mathcal{N}_{\bschi}^{\mu s}$ corepresenting these. This time the $G$-equivariant line bundles reproducing $\mu$-stability are of a particularly explicit form. With this choice $\underline{\mathcal{N}}_{\bschi}^{\mu ss}$ is corepresented by the (quasi-projective) scheme $\mathcal{N}^{\mu ss}_{\bschi}$. Moreover, the open subset $\mathcal{N}^{\mu s}_{\bschi} \subset \mathcal{N}^{\mu ss}_{\bschi}$ corepresents $\underline{\mathcal{N}}_{\bschi}^{\mu s}$ and is a coarse moduli space \cite[Thm.~4.14]{Koo1}.

\subsection{Fixed point loci of moduli spaces of sheaves}

Let $X$ be a polarized projective scheme. For any choice of Hilbert polynomial $P$, there are natural moduli functors 
\begin{align*}
\underline{\mathcal{M}}_{P}^{ss} &: (Sch/\mathbb{C})^{o} \longrightarrow Sets \\
\underline{\mathcal{M}}_{P}^{s} &: (Sch/\mathbb{C})^{o} \longrightarrow Sets 
\end{align*}
of $S$-flat families with fibres Gieseker semistable (resp.~geometrically Gieseker stable) sheaves with Hilbert polynomial $P$. See \cite[Sect.~4.1]{HL} for details. There exists a projective scheme $\mathcal{M}_{P}^{ss}$ corepresenting $\underline{\mathcal{M}}_{P}^{ss}$, an open subset $\mathcal{M}_{P}^{s} \subset \mathcal{M}_{P}^{ss}$ corepresenting $\underline{\mathcal{M}}_{P}^{s}$, and  $\mathcal{M}_{P}^{s}$ is a coarse moduli scheme \cite[Thm.~4.3.4]{HL}. Now let $X$ be a smooth projective toric variety and let $P$ have degree $\mathrm{dim}(X)$. For any $\bschi \in \X_{P}$, forgetting the $T$-equivariant structure induces a closed embedding $\mathcal{M}^{s}_{\bschi} \subset \mathcal{M}_{P}^{s}$. The action of $T$ on $X$ lifts to an action on $\mathcal{M}_{P}^{s}$ and obviously $\mathcal{M}^{s}_{\bschi} \subset \big( \mathcal{M}_{P}^{s} \big)^T$. In fact, the fixed point locus $\big( \mathcal{M}_{P}^{s} \big)^T$ can be explicitly expressed as a union of moduli spaces of $T$-equivariant sheaves.
\begin{theorem} [{{\cite[Cor.~4.10]{Koo1}}}] \label{ch. 1, sect. 4, cor. 1}
Let $X$ be a polarized smooth projective toric variety and let $P$ be a choice of Hilbert polynomial of degree $\mathrm{dim}(X)$. Then the forgetful map induces an isomorphism of schemes 
\begin{equation} \nonumber
\left( \mathcal{M}_{P}^{s} \right)^{T} \cong \coprod_{\bschi \in ( \mathcal{X}_{P} )^{\mathrm{fr}}} \mathcal{M}_{\bschi}^{s}. 
\end{equation}
\end{theorem}

\noindent Here $\left( \mathcal{X}_{P} \right)^{\mathrm{fr}} \subset \mathcal{X}_{P}$ is the collection of \emph{framed} characteristic functions. These are defined as follows. Given a $T$-equivariant torsion free sheaf $\F$ on $X$ with $\sigma$-families $\{\hat{F}^{\sigma_i}\}_{i=1, \ldots, e}$, there are unique maximally chosen integers $u_1, \ldots, u_d$ with the property 
\[
F^{\sigma_1}(\lambda_1, \ldots, \lambda_d) = 0, \ \mathrm{unless} \ \lambda_1 \geq u_1 \ \mathrm{and} \ \ldots \ \mathrm{and} \ \lambda_d \geq u_d.
\]
A characteristic function $\bschi \in \X_{P}$ is called framed if the first component $\chi^{\sigma_1}$ has the property that the integers $u_1, \ldots, u_d$ described above are all zero. For any $T$-equivariant torsion free sheaf $\F$ on $X$, there exists a \emph{unique} character $m \in M$ such that $\F \otimes \O(m)$ has framed characteristic function. Here $\O(m)$ denotes the trivial line bundle with $T$-equivariant structure induced by the character $m$. The framing ensures the forgetful map is injective. Obviously, many other choices of framing are possible.

For reflexive sheaves, there is a natural moduli functor \cite[Sect.~4.4]{Koo1}
\begin{equation*}
\underline{\mathcal{N}}_{P}^{\mu s} : (Sch/\mathbb{C})^{o} \longrightarrow Sets
\end{equation*}
of $S$-flat families with fibres geometrically $\mu$-stable reflexive sheaves with Hilbert polynomial $P$. There is an open subset $\mathcal{N}_{P}^{\mu s} \subset \mathcal{M}_{P}^{s}$ corepresenting $\underline{\mathcal{N}}_{P}^{\mu s}$ and $\mathcal{N}_{P}^{\mu s}$ is a coarse moduli space \cite[Sect.~4.4]{Koo1}. The torus action on $\mathcal{M}_{P}^{s}$ restricts to $\mathcal{N}_{P}^{\mu s}$ and the fixed point locus has the following description.
\begin{theorem} [{{\cite[Thm.~4.14]{Koo1}}}] \label{ch. 1, sect. 4, thm. 3}
Let $X$ be a polarized smooth projective toric variety and let $P$ be a choice of Hilbert polynomial of a reflexive sheaf on $X$. Then the forgetful map induces an isomorphism of schemes
\begin{equation} \nonumber
\big( \mathcal{N}_{P}^{\mu s} \big)^{T} \cong \coprod_{\bschi \in ( \mathcal{X}_{P}^{\mathrm{refl}} )^{\mathrm{fr}}} \mathcal{N}_{\bschi}^{\mu s}. 
\end{equation}
\end{theorem}

\subsection{Chern classes of equivariant sheaves}

In this paper, we want to fix the Chern classes of a sheaf rather than the Hilbert polynomial. Like in the case of Hilbert polynomial, the Chern classes of a $T$-equivariant torsion free sheaf on a toric variety are fully determined by its characteristic function. In fact, Klyachko \cite[Sect.~1.2, 1.3]{Kly4} gives an explicit formula\footnote{In the previous sections, we followed Perling's convention of ascending directions for the maps between the weight spaces as opposed to Klyachko's convention of descending directions. This results in some minus signs compared to Klyachko's original formula.}. For our purposes, we only need to know that the Chern classes are fully determined by the characteristic function, whereas the precise formula is not relevant. However, for completeness we include it.
\begin{definition} \label{ch. 1, sect. 3, def. 3}
Let $\{F(\lambda_{1}, \ldots, \lambda_{d})\}_{(\lambda_{1}, \ldots, \lambda_{d}) \in \mathbb{Z}^{d}}$ be a collection of finite-dimensional complex vector spaces. For each $i = 1, \ldots, d$, we define a $\mathbb{Z}$-linear operator $\Delta_{i}$ on the free abelian group generated by the vector spaces $\{F(\lambda_{1}, \ldots, \lambda_{d})\}_{(\lambda_{1}, \ldots, \lambda_{d}) \in \mathbb{Z}^{d}}$ determined by
\begin{equation*}
\Delta_{i}F(\lambda_{1}, \ldots, \lambda_{d}) := F(\lambda_{1}, \ldots, \lambda_{d}) - F(\lambda_{1}, \ldots, \lambda_{i-1}, \lambda_{i}-1,\lambda_{i+1}, \ldots, \lambda_{d}).
\end{equation*}
We then define $[F](\lambda_{1}, \ldots, \lambda_{d}) := \Delta_{1} \cdots \Delta_{d}F(\lambda_{1}, \ldots, \lambda_{d})$. Furthermore, we define dimension $\mathrm{dim}$ as a $\mathbb{Z}$-linear operator on the free abelian group generated by the vector spaces $\{F(\lambda_{1}, \ldots, \lambda_{d})\}_{(\lambda_{1}, \ldots, \lambda_{d}) \in \mathbb{Z}^{d}}$ in the obvious way so we can speak of $\dim [F](\lambda_{1}, \ldots, \lambda_{d})$. For example
\begin{flalign*}
\dim [F](\lambda) &= \dim F(\lambda) - \dim F(\lambda-1), \\ 
\qquad \dim [F](\lambda_{1},\lambda_{2}) &= \dim F(\lambda_{1},\lambda_{2}) - \dim F(\lambda_{1}-1,\lambda_{2}) - \dim F(\lambda_{1},\lambda_{2}-1) \\
& \ \ \ + \dim F(\lambda_{1}-1,\lambda_{2}-1). && \oslash 
\end{flalign*}  
\end{definition}
\begin{proposition}[Klyachko's Formula] \label{ch. 1, sect. 3, prop. K}
Let $X$ be a smooth projective toric variety with fan $\Delta$ and lattice $N$ of rank $d$. Let $\{\sigma_{1}, \ldots, \sigma_{e}\}$ be the cones of dimension $d$ and for each $i = 1, \ldots, e$, let $\big(\rho^{(i)}_{1}, \ldots, \rho^{(i)}_{d} \big)$ be an ordering of the rays of $\sigma_{i}$. Then any $T$-equivariant torsion free $\mathcal{F}$ on $X$ with $\sigma$-families $\{\hat{F}^{\sigma_i}\}_{i=1, \ldots, e}$ satisfies 
\begin{equation*}
\mathrm{ch}(\mathcal{F}) = \sum_{\sigma \in \Delta, \ \bslambda \in \mathbb{Z}^{\mathrm{dim}(\sigma)}} (-1)^{\mathrm{codim}(\sigma)} \dim [F^{\sigma}](\bslambda) \ \mathrm{exp}\Big( - \sum_{\rho \in \sigma(1)} \langle \bslambda, n(\rho)\rangle V(\rho) \Big).
\end{equation*} 
\end{proposition}

\noindent In this proposition, $\sigma(1)$ denotes the collection of rays of $\sigma$ and $n(\rho) \in N$ is the primitive generator of the ray $\rho$. Furthermore, $\langle \cdot, \cdot \rangle : M \times N \rightarrow \Z$ is the natural pairing and $V(\rho) \subset X$ denotes the toric divisor corresponding to the ray $\rho$. Any cone $\sigma \in \Delta$ is a face of a cone $\sigma_{i}$ of dimension $d$. Assume $\sigma$ has dimension $p$. Without loss of generality, let $\big(\rho^{(i)}_{1}, \ldots, \rho^{(i)}_{p} \big) \subset \big(\rho^{(i)}_{1}, \ldots, \rho^{(i)}_{r} \big)$ be the rays spanning $\sigma \subset \sigma_i$. Then the $\sigma$-family of the torsion free sheaf $\F|_{U_\sigma}$ is given by \cite[Prop.~2.9]{Koo1}
\[
F^{\sigma}(\lambda_{1}, \ldots, \lambda_{p}) = F^{\sigma_{i}}(\lambda_{1}, \ldots, \lambda_{p}, \infty, \ldots, \infty).
\]

\subsection{Generating functions of Euler characteristics}

In this paper, we consider the case $X = S$ is a smooth complete toric surface with polarization $H$. Instead of fixing Hilbert polynomial, we fix rank $r$ and Chern classes $c_1, c_2$. We denote by $\M_{S}^{H}(r,c_1,c_2)$ the moduli space of $\mu$-stable torsion free sheaves on $S$ with rank $r$ and Chern classes $c_1, c_2$. We want to compute the generating function of topological Euler characteristics
\[
\sum_{c_2} e(\M_{S}^{H}(r,c_1,c_2)) q^{c_2}.
\]
By the Bogomolov inequality \cite[Thm.~3.4.1]{HL}, this generating function is a formal Laurent series in $q$. Note that we compute Euler characteristics of moduli spaces of $\mu$-stable torsion free sheaves $\M_{S}^{H}(r,c_{1},c_{2})$ \emph{only} and ignore strictly $\mu$-semistables. The reason is that the descriptions of fixed point loci of Theorems \ref{ch. 1, sect. 4, cor. 1}, \ref{ch. 1, sect. 4, thm. 3} rely on simpleness in an essential way \cite{Koo1}. In the case rank and degree are coprime, i.e.~$\gcd(r,c_1 \cdot H)=1$, $\mu$-stability and Gieseker stability coincide and there are no strictly semistables, so the moduli spaces $\M_{S}^{H}(r,c_{1},c_{2})$ are projective.

For any torsion free sheaf $\F$, the natural map to its double-dual (which is reflexive \cite[Cor.~1.2]{Har}) is an injection $\F \hookrightarrow \F^{**}$ \cite[Prop.~1.1.10]{HL}. On a surface, reflexive and locally free sheaves are the same \cite[Cor.~1.4]{Har} and the cokernel of $\F \hookrightarrow \F^{**}$ is $0$-dimensional. Using this map, one can show the following \cite[Prop. 3.1]{Got3}.
\begin{proposition}  \label{GY}
Let $S$ be a smooth complete surface with polarization $H$. Let $r > 0$ and $c_{1} \in H^{2}(S,\mathbb{Z})$. Then
\begin{equation}
\sum_{c_{2}} e(\M_{S}^{H}(r, c_{1}, c_{2})) q^{c_{2}} = \frac{1}{\prod_{k=1}^{\infty}(1-q^{k})^{r e(S)}} \sum_{c_{2}} e(\N_{S}^{H}(r, c_{1}, c_{2})) q^{c_{2}}, \nonumber
\end{equation}
where $\N_{S}^{H}(r, c_{1}, c_{2})$ is the moduli space of $\mu$-stable locally free sheaves on $S$ with rank $r$ and Chern classes $c_1, c_2$.
\end{proposition} 
In the toric case, we have a torus action on the moduli spaces so $e(\N_{S}^{H}(r, c_{1}, c_{2})) = e(\N_{S}^{H}(r, c_{1}, c_{2})^T)$. Together with Theorem \ref{ch. 1, sect. 4, thm. 3}, this gives the following formula. 
\begin{proposition} \label{formula}
Let $S$ be a smooth complete toric surface with polarization $H$. Let $r > 0$ and $c_{1} \in H^{2}(S,\mathbb{Z})$. Then
\begin{equation}
\sum_{c_{2}} e(\M_{S}^{H}(r, c_{1}, c_{2})) q^{c_{2}} = \frac{1}{\prod_{k=1}^{\infty}(1-q^{k})^{r e(X)}} \sum_{c_{2}} \sum_{\bschi \in \left(\mathcal{X}_{(r,c_{1},c_{2})}^{\mathrm{refl}}\right)^{\mathrm{fr}}} e(\N_{\bschi}^{\mu s}) q^{c_{2}}, \nonumber
\end{equation}
where $\mathcal{X}_{(r,c_{1},c_{2})}^{\mathrm{refl}} \subset \X^{\mathrm{refl}}$ is the collection of characteristic function determining rank $r$ and Chern classes $c_1, c_2$ via Klyacho's formula Prop.~\ref{ch. 1, sect. 3, prop. K}.
\end{proposition}

\section{A formula for the generating function}

For any smooth complete toric surface $S$ with polarization $H$ and $r>0$, $c_1 \in H^2(S,\Z)$, we are interested in the generating function
\[
\sum_{c_2} e(\M_{S}^{H}(r,c_1,c_2)) q^{c_2}
\]
introduced in Section 2.5. In this section, we use the toric description of Proposition \ref{formula} to express this generating function in terms of Euler characteristics of certain explicit moduli spaces of stable configurations of linear subspaces (Theorem \ref{ch. 2, sect. 3, thm. 1} below). We recall that we consider $\mu$-stable torsion free sheaves \emph{only} and ignore strictly $\mu$-semistables. However, we do keep $H,r,c_1$ completely arbitrary. In the next section, we simplify the general formula of Theorem \ref{ch. 2, sect. 3, thm. 1} further in the cases: $S$ arbitrary and $r = 1$, $S = \mathbb{P}^{2}$ and $r = 1,2,3$, and $S = \mathbb{F}_{a}$ and $r = 1,2$. Here $\FF_a$ denotes the $a$th Hirzebruch surfaces and $\FF_0:= \PP^1 \times \PP^1$.

\subsection{Chern classes of equivariant locally free sheaves}

By Proposition \ref{formula}, we only need to consider reflexive, i.e.~locally free, sheaves on $S$. In this section, we compute the Chern classes of such sheaves. We start by recalling some basic facts about toric surfaces. Smooth complete toric surfaces are classified by the following proposition \cite[Sect. 2.5]{Ful}.
\begin{proposition} \label{ch. 2, sect. 3, prop. 1}
All smooth complete toric surfaces are obtained by successive blow-ups of $\mathbb{P}^{2}$ and $\mathbb{F}_{a}$ at fixed points.  
\end{proposition}

\noindent Combinatorially, such blow-ups are described by stellar subdivisions, i.e.~creating a fan $\tilde{\Delta}$ out of a fan $\Delta$ by subdividing a cone through the sum of the two primitive lattice vectors of its rays. From now on, we fix the lattice $N = \Z^2$ and let $\Delta$ be the fan of a smooth complete toric surface $S$. We denote the 2-dimensional cones by $\sigma_1, \ldots, \sigma_e$, where $e = e(S)$. We denote the rays by $\rho_1, \ldots, \rho_e$ and we let $\sigma_i$ be spanned by $\rho_i, \rho_{i+1}$. Here the index $i$ is understood modulo $e$ so $\sigma_e$ is spanned by $\rho_e, \rho_1$. Without loss of generality, we take the primitive lattice vector of $\rho_1$ to be $(1,0)$, of $\rho_2$ to be $(0,1)$, and order the rays $\rho_i$ counter-clockwise. 

The cohomology ring $H^{2 *}(S,\Z)$ can be easily described in terms of this data. First note that $H^0(S,\Z) \cong \Z$ is generated by $[S]$ and $H^4(S,\Z) \cong \Z$ by $pt$. Denote the primitive lattice vector of $\rho_i$ by $n_i$ and denote the toric divisor corresponding to $\rho_i$ by $D_i$. Then $H^2(S,\Z)$ is generated by $D_1, \ldots, D_e$ modulo the relations \cite[Sect. 5.2]{Ful}
\begin{align*}
&D_{1} + \sum_{i=3}^{e} \langle (1,0), n_{i} \rangle D_{i} = 0, \\ 
&D_{2} + \sum_{i=3}^{e} \langle (0,1), n_{i} \rangle D_{i} = 0.
\end{align*}
Here $M = \Z^2$ and $\langle \cdot, \cdot \rangle$ is the standard inner product. By \cite[Sect.~2.5]{Ful}, $D_i D_j = 0$ unless $j=i+1$ and 
\[
D_{1}D_{2} = D_{2}D_{3} = \cdots = D_{e-1}D_{e} = D_{e}D_{1} = pt.
\]
Finally, the self-intersections $D_{i}^{2} = -a_{i}$ are determined by the equation $n_{i-1} + n_{i+1} = a_{i}n_{i}$ \cite[Sect. 2.5]{Ful}. For future reference, it is convenient to define $\xi_{i} := - \langle (1,0), n_{i} \rangle$ and $\eta_{i} := - \langle (0,1), n_{i} \rangle$. Note that the integers $\{a_{i}\}_{i=1}^{e}$, $\{\xi_{i}\}_{i=3}^{e}$, $\{\eta_{i}\}_{i=3}^{e}$ are entirely determined by the fan $\Delta$.

By Theorem \ref{Kly}, a $T$-equivariant rank $r$ torsion free sheaf $\F$ on $S$ is described by multifiltrations $\{F^{\sigma_{i}}(\lambda_{1},\lambda_{2})\}_{i=1, \ldots, e}$ of $\C^{\oplus r}$ satisfying the gluing conditions
\begin{align} \label{glue1}
F^{\sigma_{i}}(\infty, \lambda) = F^{\sigma_{i+1}}(\lambda, \infty), \ \mathrm{for \ all \ } \lambda \in \mathbb{Z}. 
\end{align}
Moreover, a $T$-equivariant rank $r$ locally free sheaf $\F$ on $S$ is simply described by flags $\{V^{\rho_i}(\lambda)\}_{i = 1, \ldots, e}$ of $\C^{\oplus r}$ (Section 2.1). As we discussed, the corresponding $\sigma$-families are defined by
\begin{equation*} 
F^{\sigma_i}(\lambda_1, \lambda_2) := V^{\rho_i}(\lambda_1) \cap V^{\rho_{i+1}}(\lambda_2).
\end{equation*}
The flags $\{V^{\rho_i}(\lambda)\}_{i = 1, \ldots, e}$ can be described by indicating the integers where the vector spaces jump together with the subspaces occurring in the flag. More precisely, for each $i = 1, \ldots, e$, there exist unique integers $u_i \in \Z$, $v_{1,i}, \ldots, v_{r-1,i} \in \Z_{\geq 0}$ and subspaces $p_{1,i} \in \mathrm{Gr}(1,r), \ldots, p_{r-1,i} \in \mathrm{Gr}(r-1,r)$ such that
\begin{equation}
\begin{split} \label{toricdataformula}
V^{\rho_{i}}(\lambda) = \left\{\begin{array}{cc}  0 & \mathrm{if \ } \lambda < u_{i} \\ p_{1,i} & \mathrm{if \ } u_{i} \leq \lambda < u_{i} + v_{1,i} \\ p_{2,i} & \mathrm{if \ } u_{i} + v_{1,i} \leq \lambda < u_{i} + v_{1,i} + v_{2,i} \\ \ldots & \ldots \\ \mathbb{C}^{\oplus r} & \mathrm{if \ } u_{i} + v_{1,i} +  \ldots + v_{r-1,i} \leq \lambda. \end{array} \right. 
\end{split}
\end{equation}
Note that $v_{a,i}$ could be zero in which case $p_{a,i}$ does not occur. At such places, the flag jumps more than 1 in dimension. 

\begin{definition} \label{toricdata}
Instead of describing a $T$-equivariant locally free sheaf $\F$ on $S$ by the flags $\{V^{\rho_i}(\lambda)\}_{i = 1, \ldots, e}$, we can also describe it by the data $\{(u_i, v_{a,i}, p_{a,i})\}_{a=1, \ldots, r-1, i=1, \ldots, e}$ introduced above. We refer to $\{(u_i, v_{a,i}, p_{a,i})\}_{a=1, \ldots, r-1, i=1, \ldots, e}$ as \emph{toric data} and abbreviate it by $({\bf{u}}, {\bf{v}}, {\bf{p}})$. \hfill $\oslash$ 
\end{definition}
\begin{proposition} \label{ch. 2, sect. 3, prop. 2}
Let $\mathcal{F}$ be a $T$-equivariant rank $r$ locally free sheaf on $S$ described by toric data $({\bf{u}}, {\bf{v}}, {\bf{p}})$. Then
\begin{align}
c_1(\mathcal{F}) =&- \sum_{i=1}^{e}\big( r u_{i} + \sum_{a=1}^{r-1}(r-a) v_{a,i} \big)D_{i}, \nonumber \\
\mathrm{ch}_2(\mathcal{F}) =& \frac{1}{2} \big(\sum_{i=1}^{e} u_{i}D_{i} \big)^{2} + \frac{1}{2} \sum_{a=1}^{r-1} \Big(\sum_{i=1}^{e}\big(u_{i} + \sum_{b=1}^{a} v_{b,i} \big)D_{i} \Big)^{2} \nonumber \\
&- \sum_{i=1}^{e} \sum_{a,b=1}^{r-1} v_{a,i} v_{b,i+1} \big( \min\{a,b\}-\dim(p_{a,i} \cap p_{b,i+1}) \big) \ pt. \nonumber
\end{align}
\end{proposition}
\begin{proof}
In the case $r=1$, the sheaf $\F$ is a line bundle and described by integers $u_1, \ldots, u_e$ (Section 2.1). It is easy to see that \cite[Sect.~4.2]{Koo1} 
\[
c_1(\F) = -\sum_{i=1}^e u_i D_i.
\]
Therefore
\begin{align}
\mathrm{ch}(\mathcal{F}) = \mathrm{exp} \Big(-\sum_{i=1}^{e} u_{i} D_{i} \Big) =1 - \sum_{i=1}^{e} u_{i} D_{i} + \frac{1}{2} \big( \sum_{i=1}^{e} u_{i} D_{i} \big)^{2}. \nonumber
\end{align}

In the case $r>0$ and $p_{a,i} = p_{a,i+1}$ for all $a,i$, the sheaf $\F$ is a direct sum of $T$-equivariant line bundles 
\[
\F = \bigoplus_{a=1}^{r} \L_a.
\]
Here $\L_a$ is defined by flags $\{L_{a}^{\rho_i}(\lambda)\}_{i=1, \ldots, e}$, where $L_{a}^{\rho_i}(\lambda) = \C$ if $\lambda \geq u_i + \sum_{b=1}^{a-1} v_{b,i}$ and $L_{a}^{\rho_i}(\lambda) = 0$ otherwise. This immediately implies the following formula  
\begin{align}
&\ch(\F) = \sum_a \ch(\L_a) \label{directsum} \\
&= r - \sum_{i=1}^{e}\big( r u_{i} + \sum_{a=1}^{r-1}(r-a) v_{a,i} \big)D_{i} + \frac{1}{2} \big(\sum_{i=1}^{e} u_{i}D_{i} \big)^{2} + \frac{1}{2} \sum_{a=1}^{r-1} \Big(\sum_{i=1}^{e}\big(u_{i} + \sum_{b=1}^{a} v_{b,i} \big)D_{i} \Big)^{2}. \nonumber
\end{align}

For the general case, we use Klyachko's formula (Proposition \ref{ch. 1, sect. 3, prop. K}). Actually, we do not need the precise form of the formula, but merely observe $\ch(\F)$ only depends on the characteristic function $\bschi_{\F}$ (Definition \ref{ch. 1, sect. 3, def. 1}). For each $a = 1, \ldots, r$, define a $T$-equivariant torsion free subsheaf $\G_a \subset \L_a$ by the following $\sigma$-families $\{G^{\sigma_i}_{a}(\lambda_1, \lambda_2)\}_{i=1, \ldots, e}$
\begin{equation}
G_{a}^{\sigma_{i}}(\lambda_{1}, \lambda_{2}) = \left\{ \begin{array}{cc} \mathbb{C} & \mathrm{if} \ \mathrm{dim}(F^{\sigma_{i}}(\lambda_{1}, \lambda_{2})) \geq a \\ 0 & \mathrm{otherwise.} \end{array} \right. \nonumber
\end{equation}
Then by construction $\bschi_{\F} = \sum_a \bschi_{\G_a} = \bschi_{\bigoplus_a \G_a}$ so $\ch(\F) = \ch(\bigoplus_a \G_a)$. The sheaf $\bigoplus_a \G_a$ is a $T$-equivariant subsheaf of $\bigoplus_a \L_a$ with 0-dimensional cokernel $\cQ$. The length of $\cQ$ is easily seen to be 
\[
\sum_{i=1}^{e} \sum_{a,b=1}^{r-1} v_{a,i} v_{b,i+1} \big( \min\{a,b\}-\dim(p_{a,i} \cap p_{b,i+1}) \big).
\]
Subtracting this from equation (\ref{directsum}) gives the answer.
\end{proof}

\subsection{Main theorem}

\noindent \emph{Characteristic functions of locally free sheaves.} By Proposition \ref{ch. 2, sect. 3, prop. 2}, we now know how a characteristic function $\bschi \in \X^{\mathrm{refl}}$ determines rank and Chern classes. Next, we want to say a bit more about $\bschi$ itself. Let $\sigma_i \in \Delta$ be a cone of maximal dimension and consider the corresponding $T$-invariant affine open subset $U_{\sigma_i}$. Let $\F$ be a $T$-equivariant locally free sheaf of rank $r$ on $S$. The restriction $\F|_{U_{\sigma_i}}$ splits into a sum of $T$-equivariant line bundles on $U_{\sigma_i}$
\[
\F|_{U_{\sigma_i}} \cong \bigoplus_{a=1}^r \L_a.
\]
Note that in general, we do not have such a splitting globally. From this splitting, we can read off the $i$th component $\chi^{\sigma_i}_{\F}$ of the characteristic function $\bschi_{\F}$. Indeed, let $\L_a$ be generated by a homogeneous element with character $m^{\sigma_i}_{a}$, then the collection of characters $\{m^{\sigma_i}_1, \ldots, m^{\sigma_i}_r\}$ completely determine $\chi^{\sigma_i}_{\F}$. Let us make this explicit. As before, denote the primitive generator of ray $\rho_i$ by $n_i$ and the pairing by $\langle \cdot, \cdot \rangle$. Define the Heaviside function
\begin{align*}
&H_{m^{\sigma_i}_{a}} : M \rightarrow \Z, \\
&H_{m^{\sigma_i}_{a}} (\lambda_1,\lambda_2) = \left\{ \begin{array}{cc} 1 & \mathrm{if} \ \lambda_1 \geq \langle m^{\sigma_i}_{a}, n_i \rangle \ \mathrm{and} \ \lambda_2 \geq \langle m^{\sigma_i}_{a}, n_{i+1} \rangle \\ 0 & \mathrm{otherwise.}  \end{array}\right.
\end{align*}
Recall that we use the primitive generators $(n_i,n_{i+1})$ as a basis for $N$ and the dual basis as a basis for $M$ (Sections 2.1 and 3.1). Then
\[
\chi^{\sigma_i}_{\F} = \sum_{a=1}^r H_{m^{\sigma_i}_{a}}.
\]
So indeed $\{m^{\sigma_i}_1, \ldots, m^{\sigma_i}_r\}$ fully determines $\chi^{\sigma_i}$ and vice versa. By the gluing conditions (\ref{glue1}), a sequence $\{\{m^{\sigma_i}_1, \ldots, m^{\sigma_i}_r\}\}_{i=1, \ldots, e}$ determines a characteristic function of a rank $r$ $T$-equivariant locally free sheaf on $S$ if and only if
\begin{equation*} 
\langle m^{\sigma_i}_{a}, n_{i+1} \rangle = \langle m^{\sigma_{i+1}}_{a}, n_{i+1} \rangle, 
\end{equation*}
for all $a=1, \ldots, r$ and $i=1, \ldots, e$. 

Now let $\F$ be any $T$-equivariant locally free sheaf on $S$ with characteristic function $\{\{m^{\sigma_i}_1, \ldots, m^{\sigma_i}_r\}\}_i$ and toric data $({\bf{u}}, {\bf{v}}, {\bf{p}})$. The notion of toric data was introduced in Definition \ref{toricdata}. The integers $u_i$, $v_{a,i}$ are full determined by the characteristic function via the following equations
\begin{align*}
\langle m^{\sigma_i}_{1}, n_i \rangle &= \langle m^{\sigma_{i-1}}_{1}, n_{i} \rangle = u_i, \\
\langle m^{\sigma_i}_{2}, n_i \rangle &= \langle m^{\sigma_{i-1}}_{2}, n_i \rangle = u_i + v_{1,i}, \\
& \cdots \\
\langle m^{\sigma_i}_{r}, n_i \rangle &= \langle m^{\sigma_{i-1}}_{r}, n_i \rangle = u_i + v_{1,i} + \cdots + v_{r-1,i}.
\end{align*}
Although a characteristic function does not determine the continuous parameters $p_{a,i}$, it does determine the dimensions
\[
\dim( p_{a,i} \cap p_{b,i+1} ),
\]
for all $i=1, \ldots e$ and $a,b=1, \ldots, r-1$. We denote these dimensions by $\delta_{a,b,i}:= \dim( p_{a,i} \cap p_{b,i+1} )$. Note that $\delta_{a,b,i} \in \{0,1, \ldots, \min\{a,b\}\}$. We abbreviate the data $\{(u_i,v_{a,i},\delta_{a,b,i})\}_{a,b,i}$ by $({\bf{u}}, {\bf{v}}, \bsdelta)$. Clearly the data of a characteristic function $\bschi$ is equivalent to the data $({\bf{u}}, {\bf{v}}, \bsdelta)$. From now on, we identify the two notions 
\[
\bschi \leftrightarrow ({\bf{u}}, {\bf{v}}, \bsdelta).
\]
The reason for introducing this notation is because Proposition \ref{ch. 2, sect. 3, prop. 2} expresses the Chern classes of a $T$-equivariant locally free sheaf with characteristic function $({\bf{u}}, {\bf{v}}, \bsdelta)$ in terms of this data. \\

\noindent \emph{Stratification.} As we have seen in Sections 2.1 and 3.1, $T$-equivariant locally free sheaves of rank $r$ on $S$ are described by toric data $({\bf{u}},{\bf{v}},{\bf{p}})$. Such toric data is naturally parametrized by the closed points of the following variety
\begin{equation} \label{flag}
\coprod_{u_1, \ldots, u_e \in \Z} \coprod_{{\footnotesize{\begin{array}{c} v_{1,1}, \ldots, v_{r-1,1} \geq 0 \\ \ldots \\ v_{1,e}, \ldots, v_{r-1,e} \geq 0 \end{array}}}} \prod_{i=1}^e \mathrm{Flag}(u_i, v_{1,i}, \ldots, v_{r-1,i}),
\end{equation}
where $\mathrm{Flag}(u_i, v_{1,i}, \ldots, v_{r-1,i})$ is the partial flag variety of flags $p_{1,i} \subset \cdots \subset p_{r-1,i} \subset \C^{\oplus r}$. The labels $u_i$, $v_{a,i}$ allow us to recover the toric data by formula (\ref{toricdataformula}). For any $({\bf{u}}, {\bf{v}}, \bsdelta) \in \X^{\mathrm{refl}}$, we denote by $\D_{({\bf{u}}, {\bf{v}}, \bsdelta)}$ the collection of toric data $({\bf{u}}, {\bf{v}}, {\bf{p}})$ with characteristic function $({\bf{u}}, {\bf{v}}, \bsdelta)$. Clearly, $\D_{({\bf{u}}, {\bf{v}}, \bsdelta)}$ is naturally a locally closed\footnote{Note that for any finite product of Grassmannians $\prod_{i} \mathrm{Gr}(n_{i},N)$, the map $\{p_{i}\}_i \mapsto \mathrm{dim}\left( \bigcap_{i} p_{i} \right)$ is upper semicontinuous.} subset of (\ref{flag}). We can now stratify (\ref{flag}) as follows
\begin{equation*} 
\coprod_{u_1, \ldots, u_e \in \Z} \coprod_{{\footnotesize{\begin{array}{c} v_{1,1}, \ldots, v_{r-1,1} \geq 0 \\ \ldots \\ v_{1,e}, \ldots, v_{r-1,e} \geq 0 \end{array}}}} \coprod_{{\footnotesize{\begin{array}{c} \delta_{a,b,i} \in \{0,1, \ldots, \min\{a,b\} \} \\ \mathrm{for \ all} \ i=1, \ldots, e \\ \mathrm{and} \ a,b=1, \ldots, r-1 \end{array}}}} \D_{({\bf{u}}, {\bf{v}}, \bsdelta)}.
\end{equation*}
The advantage of this stratification is that any $T$-equivariant locally free sheaf on $S$ with toric data in $\D_{({\bf{u}}, {\bf{v}}, \bsdelta)}$ has the \emph{same} Chern character by Proposition \ref{ch. 2, sect. 3, prop. 2}.

Each component of the variety (\ref{flag}) is naturally a closed subscheme of 
\[
\prod_{i=1}^{e} \prod_{a=1}^{r-1} \mathrm{Gr}(a,r),
\]
where we omit the factor indexed by $a$, $i$ when $v_{a,i} = 0$. This product of Grassmannians carries a natural action of $\mathrm{SL}(r,\C)$, which keeps each factor $\D_{({\bf{u}}, {\bf{v}}, \bsdelta)}$ invariant. Equivariant isomorphism classes of ample linearizations on $\prod_{i=1}^{e} \prod_{a=1}^{r-1} \mathrm{Gr}(a,r)$ are in 1-1 correspondence with sequences of positive integers $\{\kappa_{a,i}\}_{a = 1, \ldots, r-1, i=1, \ldots, e}$ by \cite[Sect.~11.1]{Dol}. On a factor $\D_{({\bf{u}}, {\bf{v}}, \bsdelta)}$, we are interested in the following linearization. The toric data in $\D_{({\bf{u}}, {\bf{v}}, \bsdelta)}$ gives rise to integers $u_i$, $v_{a,i}$ and we take the ample linearization
\[
\{(H \cdot D_{i}) v_{a,i} \}_{a,i}
\]
on the product of Grassmannians and restrict it to $\D_{({\bf{u}}, {\bf{v}}, \bsdelta)}$. Recall that $H$ is the (fixed) polarization on $S$ and the $D_i$ are the toric divisors (Section 3.1). It is proved in \cite[Prop.~3.20]{Koo1}, that the notion of GIT stability on $\D_{({\bf{u}}, {\bf{v}}, \bsdelta)}$ we obtain in this way coincides with $\mu$-stability. More precisely, any $T$-equivariant locally free sheaf $\F$ on $S$ with toric data in $({\bf{u}}, {\bf{v}}, {\bf{p}}) \in \D_{({\bf{u}}, {\bf{v}}, \bsdelta)}$ is $\mu$-semistable if and only if $({\bf{u}}, {\bf{v}}, {\bf{p}})$ corresponds to a GIT semistable point and $\F$ is $\mu$-stable if and only if $({\bf{u}}, {\bf{v}}, {\bf{p}})$ corresponds to a properly GIT stable point (with respect to the chosen linearization). The previous discussion combined with Theorem \ref{ch. 1, sect. 4, thm. 3} gives the following proposition. 
\begin{proposition} \label{ch. 2, sect. 3, prop. 5}
Let $S$ be a smooth complete toric surface with polarizarion $H$. Let $r > 0$ and $c_{1} \in H^{2}(S,\mathbb{Z})$. Then for any $c_{2} \in H^4(S,\Z) \cong \mathbb{Z}$, there is a canonical isomorphism
\begin{align*}
\N_{S}^{H}(r,c_{1},c_{2})^{T} \cong \coprod_{{\footnotesize{\begin{array}{c} u_i, v_{a,i} \\ \mathrm{giving \ rise \ to \ } c_{1} \end{array}}}} \coprod_{{\footnotesize{\begin{array}{c} \delta_{a,b,i} \\ \mathrm{giving \ rise \ to \ } c_{2} \end{array}}}} \mathcal{D}_{({\bf{u}}, {\bf{v}}, \bsdelta)}^{s} / \mathrm{SL}(r,\mathbb{C}), \nonumber
\end{align*}
where $\mathcal{D}_{({\bf{u}}, {\bf{v}}, \bsdelta)}^{s} \subset \mathcal{D}_{({\bf{u}}, {\bf{v}}, \bsdelta)}$ is the open subset of properly GIT stable points with respect to the polarization $\{(H \cdot D_{i}) v_{a,i} \}_{a,i}$ and the quotients are good geometric quotients. 
\end{proposition}

\noindent Some comments about this proposition are in order. Firstly, in the union over $u_1, \ldots, u_e \in \Z$ we take $u_1=u_2=0$ and $u_3, \ldots, u_e \in \Z$ arbitrary. This is because the disjoint union in Theorem \ref{ch. 1, sect. 4, thm. 3} is over \emph{framed} characteristic functions. Secondly, we note that it makes sense to speak of $u_i$, $v_{a,i}$ giving rise to some fixed $c_1 \in H^2(S,\Z)$ by the formula of Proposition \ref{ch. 2, sect. 3, prop. 2}.  Thirdly, by the same proposition, it makes sense to speak of $u_i$, $v_{a,i}$, $\delta_{a,b,i}$ giving rise to some fixed $c_2 \in H^4(S,\Z) \cong \Z$. \\

\noindent \emph{Main theorem.} We introduce some final notation. For a fixed $c_{1} = \sum_{i=3}^{e} f_{i} D_{i} \in H^{2}(X,\mathbb{Z})$, we define
\[
C:= \big\{ \{v_{a,i}\}_{a,i} \in \Z_{\geq 0}^{(r-1)e} \ : \  r \ | \ -f_{i} + \sum_{a=1}^{r-1} a \left(v_{a,1}\xi_{i} + v_{a,2}\eta_{i} + v_{a,i} \right) \ \forall \ i=3, \ldots e \big\}.
\]
We suppress the dependence of $C$ on $S$, $r$, and $c_1$ as we think of these as fixed. Recall that the integers $\xi_i$, $\eta_i$ were introduced in Section 3.1 and are entirely determined by the fan of $S$. We also introduce the following complicated quadratic polynomial in the variables ${\bf{v}}=\{v_{a,i}\}_{a,i}$
\begin{align*}
&Q({\bf{v}}):= \\
&\frac{1}{2} \Big( \sum_{i=3}^{e} f_{i}D_{i} \Big)^{2} - \frac{1}{2r^{2}} \sum_{a=0}^{r-1} \Bigg[ \sum_{i=3}^{e} \Big( -f_{i} - \sum_{b=1}^{r-1}(r-b) v_{b,i} + \big\{-\sum_{b=1}^{r-1} (r-b) v_{b,1} + \sum_{b=1}^{a} r v_{b,1} \big\} \xi_{i} \\ 
&+ \big\{-\sum_{b=1}^{r-1} (r-b) v_{b,2} + \sum_{b=1}^{a} r v_{b,2} \big\} \eta_{i} + \sum_{b=1}^{a} r v_{b,i} \Big)D_{i} \Bigg]^{2}. 
\end{align*}
As before, we suppress the dependence of $Q$ on $S$, $r$, $c_1$. For any ${\bf{v}} = \{v_{a,i}\}_{a,i} \in C$ and $u_1=u_2=0$, there are \emph{unique} $u_3, \ldots, u_e$ such that $u_i$, $v_{a,i}$ determine $c_1$ by the formula of Proposition \ref{ch. 2, sect. 3, prop. 2}. For any choice of $\bsdelta = \{\delta_{a,b,i}\}_{a,b,i}$ we define
\begin{align*}
R({\bf{v}}, \bsdelta) &:= \sum_{i=1}^e \sum_{a,b=1}^{r-1} v_{a,i} v_{b,i+1} \big( \min\{a,b\}-\delta_{a,b,i} \big), \\
\D_{({\bf{v}},\bsdelta)} &:= \D_{({\bf{u}},{\bf{v}},\bsdelta)}, \ \mathrm{where} \ u_1=u_2=0 \ \mathrm{and} \ u_3, \ldots, u_e \ \mathrm{determine} \ c_1.
\end{align*}
Combining Propositions \ref{formula}, \ref{ch. 2, sect. 3, prop. 2}, \ref{ch. 2, sect. 3, prop. 5} gives the following explicit formula for the generating function for any $S$, $H$, $r$, and $c_1$. 
\begin{theorem} \label{ch. 2, sect. 3, thm. 1}
Let $S$ be a smooth complete toric surface with polarizarion $H$. Let $r >0$ and $c_{1} = \sum_{i=3}^{e} f_{i} D_{i} \in H^{2}(S,\mathbb{Z})$. Then 
\[
\sum_{c_{2}} e(\M_{S}^{H}(r,c_{1},c_{2})) q^{c_{2}} = \frac{1}{\prod_{k=1}^{\infty}(1-q^{k})^{r e(X)}} \sum_{{\bf{v}} \in C} \sum_{\bsdelta} e(\mathcal{D}_{({\bf{v}},\bsdelta)}^{s} / \mathrm{SL}(r,\mathbb{C})) \ q^{Q({\bf{v}}) + R({\bf{v}}, \bsdelta)},
\]
where $\mathcal{D}_{({\bf{v}},\bsdelta)}^{s} \subset \D_{({\bf{v}},\bsdelta)}$ is the open subset of properly GIT stable point with respect to the polarization $\{(H \cdot D_{i}) v_{a,i} \}_{a,i}$ and the quotients are good geometric quotient.
\end{theorem}

\section{Examples}

In this section we specialize the expression of Theorem \ref{ch. 2, sect. 3, thm. 1} to the following cases: any $S$ and $r=1$, $S=\PP^2$ and $r=1,2,3$, and $S = \mathbb{F}_a$ and $r=2$. Some of these cases have been considered individually by various authors including Ellingsrud and Str{\o}mme, G\"ottsche, Klyachko, Yoshioka and Weist. In the case $S = \mathbb{F}_a$, we study the dependence on the choice of polarization and compare to Joyce's general theory of wall-crossing for motivic invariants counting (semi)stable objects in an abelian category.

The case of any toric surface $S$ and $r=1$ trivially gives 
\begin{equation} 
\sum_{c_{2}} e(\M_{S}(1,c_{1},c_{2})) q^{c_{2}} =  \frac{1}{\prod_{k=1}^{\infty}(1-q^{k})^{e(X)}}. \nonumber
\end{equation}
For any (not necessarily toric) surface $S$, we have $\M_{S}(1,c_{1},c_{2}) \cong  \mathrm{Pic}^0(S) \times \Hilb^{c_2}(S)$, where $\Hilb^{c_2}(S)$ is the Hilbert scheme of $c_2$ points on $S$ and $\Pic^0(S)$ is the Picard torus of $S$. Therefore, the above is also the generating function of Euler characteristics of Hilbert schemes of points on $S$. For $S = \PP^2$ or $\mathbb{F}_a$, Ellingsrud and Str{\o}mme \cite{ES} computed the Betti numbers of $\mathrm{Hilb}^n(S)$ using localization techniques. Subsequently, G\"ottsche \cite{Got1} computed the Betti numbers of $\mathrm{Hilb}^n(S)$ for any smooth complete surface $S$. His proof uses the Weil conjectures.

\subsection{Rank 2 on $\mathbb{P}^{2}$ and $\mathbb{F}_{a}$}

In the $r=2$ case, the expression of Theorem \ref{ch. 2, sect. 3, thm. 1} involves Euler characteristics of configuration spaces of points on $\PP^1$. Note that these configuration spaces depend explicitly on the choice of polarization $H$ on $S$. For the toric data (Definition \ref{toricdata}) $({\bf{u}},{\bf{v}},{\bf{p}})$ of a rank 2 locally free sheaf on $S$, we define $v_i:=v_{1,i}$ and $p_i:=p_{1,i}$. For the characteristic function (see Section 3.2) $({\bf{u}},{\bf{v}},\bsdelta)$ of such a sheaf, we moreover write $\delta_i:=\delta_{1,i}$.

\subsubsection{Rank 2 on $\mathbb{P}^{2}$}

Let $S = \PP^2$. The generating function does not depend on choice of polarization, so we suppress it from the notation. Since $e(S)=3$ and $r=2$, the spaces $\D_{({\bf{v}},\bsdelta)}$ of Theorem \ref{ch. 2, sect. 3, thm. 1} are locally closed subsets of $(\PP^1)^3$. The only possibly non-empty quotients $\D_{({\bf{v}},\bsdelta)}^{s} / \mathrm{SL}(2,\C)$ are those for which all $v_i>0$ and all $\delta_{i}$ are 0. In this case 
\[
\D_{({\bf{v}},{\bf{0}})}^{s} \subset (\PP^1)^3
\]
is the open subset of triples $(p_1,p_2,p_3)$ with all $p_i$ mutually distinct. The quotients $\D_{({\bf{v}},{\bf{0}})}^{s} / \mathrm{SL}(2,\C)$ are either empty or consist of one reduced point depending on the value of the polarization. Specifically
\begin{equation*}
\D_{({\bf{v}},{\bf{0}})}^{s} / \mathrm{SL}(2,\C) = \left\{\begin{array}{cc} pt & \mathrm{if \ } v_i < v_j + v_k \ \mathrm{for \ all \ } \{i,j,k\}=\{1,2,3\} \\ 0 & \mathrm{otherwise.} \end{array} \right. 
\end{equation*} 
The notation ``for all $\{i,j,k\}=\{1,2,3\}$'' means ``for all $i \in \{1,2,3\}$, $j \in \{1,2,3\} \setminus \{i\}$, and $k \in \{1,2,3\} \setminus \{i,j\}$''. Writing the first Chern class as $c_1 = f H$, where $H$ is the hyperplane class, Theorem \ref{ch. 2, sect. 3, thm. 1} gives
\begin{align}
&\prod_{k=1}^{\infty} (1-q^{k})^{6} \sum_{c_{2}} e(\M_{\PP^2}(2,c_{1},c_{2})) q^{c_{2}} = \sum_{{\footnotesize{\begin{array}{c} v_{1}, v_{2}, v_{3} > 0 \ \mathrm{s.t.} \\ 2 \ | \ -f + v_{1} + v_{2} + v_{3} \\ v_{1} < v_{2} + v_{3} \\ v_{2} < v_{1} + v_{3} \\ v_{3} < v_{1} + v_{2} \end{array}}}} q^{\frac{f^{2}}{4} + \frac{1}{2} \sum_{i<j} v_i v_j - \frac{1}{4} \sum_i v_{i}^{2}}. \label{eqnp2}
\end{align}

Let $S$ be any smooth complete surface, $H$ an ample divisor, $r > 0$, $c_{1} \in H^{2}(S,\mathbb{Z})$ and $c_{2} \in H^4(S,\Z) \cong \mathbb{Z}$. Let $a$ be a Weil divisor. Applying $- \otimes \mathcal{O}_{S}(a)$, we obtain an isomorphism
\begin{equation}
\M_{S}^{H}(r,c_{1},c_{2}) \cong \M_{S}^{H}(r,c_{1}+ra,(r-1)c_{1}a+\frac{1}{2}r(r-1)a^{2}+c_{2}). \nonumber
\end{equation}
This uses the fact that $- \otimes \mathcal{O}_{S}(a)$ preserves $\mu$-stability. We deduce
\begin{equation}
\sum_{c_{2}} e(\M_{S}^{H}(r,c_{1}+ra,c_{2}))q^{c_{2}} = q^{(r-1)c_{1}a+\frac{1}{2}r(r-1)a^{2}} \sum_{c_{2}} e(\M_{S}^{H}(r,c_{1},c_{2}))q^{c_{2}}. \label{ch. 2, eqn2}
\end{equation}
So for $S = \mathbb{P}^{2}$ and $r = 2$, the only two interesting values for $c_1$ are $0$ and $1$.  
\begin{corollary} \label{ch. 2, sect. 4, cor. 2}
On $S = \PP^2$, we have the following rank 2 generating functions for Euler characteristics of moduli spaces of $\mu$-stable torsion free sheaves
\begin{align}
&\sum_{c_{2}} e(\M_{\PP^2}(2,0,c_{2})) q^{c_{2}} = \frac{1}{\prod_{k=1}^{\infty}(1-q^{k})^{6}} \sum_{m = 1}^{\infty} \sum_{n = 1}^{\infty} \frac{q^{mn+m+n}}{1-q^{m+n}}, \nonumber \\
&\sum_{c_{2}} e(\M_{\PP^2}(2,1,c_{2})) q^{c_{2}} = \frac{1}{\prod_{k=1}^{\infty}(1-q^{k})^{6}} \sum_{m = 1}^{\infty} \sum_{n = 1}^{\infty} \frac{q^{mn}}{1-q^{m+n-1}}. \nonumber 
\end{align}
\end{corollary}
\begin{proof}
The corollary follows from rewriting equation (\ref{eqnp2}). Using the substitutions $\xi = \frac{1}{2}(v_{1}+v_{2}-v_{3})$, $\eta = \frac{1}{2}(v_{1}-v_{2}+v_{3})$, $\zeta = \frac{1}{2}(-v_{1}+v_{2}+v_{3})$, the set
\begin{align}
&\big\{ (v_{1}, v_{2}, v_{3}) \in \mathbb{Z}^3 \ : \ 2 \ | \ -f+v_{1}+v_{2}+v_{3}, \ v_{i}>0, \ v_{i} < v_{j}+v_{k} \ \forall \ \{i,j,k\} = \{1,2,3\} \big\} \nonumber 
\end{align}
becomes
\begin{equation}
\big\{ (\xi, \eta, \zeta) \in \mathbb{Q}_{>0}^{3} \ : \ 2 \ | \ -f + 2\xi+2\eta+2\zeta, \ \xi+\eta \in \mathbb{Z}, \ \xi + \zeta \in \mathbb{Z}, \ \eta + \zeta \in \mathbb{Z} \big\}. \nonumber
\end{equation}
Using the substitutions $\xi = \frac{2k-f}{2}$, $\eta = m - \frac{2k-f}{2}$, $\zeta = n - \frac{2k-f}{2}$, this set becomes
\begin{equation}
\big\{ (k,m,n) \in \mathbb{Z}^3 \ : \ k > \frac{f}{2}, \ m > k - \frac{f}{2}, \ n > k - \frac{f}{2} \big\}. \nonumber
\end{equation}
Applying these substitutions and setting $f=1$ gives  
\begin{equation*}
\sum_{c_{2}} e(\M_{\PP^2}(2,1,c_{2})) q^{c_{2}} = \frac{1}{\prod_{p=1}^{\infty}(1-q^{p})^{6}} \sum_{k=1}^{\infty} \sum_{m=k}^{\infty} \sum_{n=k}^{\infty} q^{mn-k(k-1)}, 
\end{equation*}
and a similar formula holds for $c_1=0$. The result follows from the geometric series.
\end{proof}

\noindent \emph{Comparison to existing literature.} In \cite{Yos}, Yoshioka derives an expression for the generating function of Poincar\'e polynomials of $\M_{\PP^2}(2,1,c_{2})$ using the Weil Conjectures. Specializing his formula to Euler characteristics gives
\begin{align*}
&\sum_{c_{2}} e(\M_{\PP^2}(2,1,c_{2})) q^{c_{2}} = \\
&\frac{1}{\prod_{k=1}^{\infty}(1-q^{k})^{6}} \left(\frac{1}{2\sum_{m \in \mathbb{Z}} q^{m^{2}}}\right) \sum_{n=0}^{\infty}\left(\frac{2-4n}{1-q^{2n+1}} + \frac{8q^{2n+1}}{(1-q^{2n+1})^{2}} \right)q^{(n+1)^{2}}. 
\end{align*}
Equating to the formula of Corollary \ref{ch. 2, sect. 4, cor. 2} gives an interesting identity of formal power series. Although it does not seem to be easy to show the equality directly, one can numerically check agreement of the coefficients up to large order. 

In \cite{Kly4}, Klyachko computes $\sum_{c_{2}} e(\M_{\PP^2}(2,1,c_{2}))q^{c_{2}}$ and our paper basically follows his philosophy. In fact, the prequel to this paper \cite{Koo1} lays the foundations of many ideas appearing in \cite{Kly4} in the case of pure sheaves of any dimension on any smooth toric variety. This paper can be seen as a systematic application of these ideas to smooth toric surfaces. Klyachko expresses his answer as
\begin{equation}
\sum_{c_{2}} e(\M_{\PP^2}(2,1,c_{2})) q^{c_{2}} = \frac{1}{\prod_{k=1}^{\infty}(1-q^{k})^{6}} \sum_{m=1}^{\infty} 3 H(4m-1) q^{m}, \nonumber
\end{equation}
where $H(D)$ is the Hurwitz class number
\begin{equation}
H(D) = \left( \begin{array}{c} \mathrm{number \ of \ integer \ binary \ quadratic \ forms \ } Q \mathrm{ \ of} \\ \mathrm{discriminant} -D \mathrm{ \ counted \ with \ weight \ } \frac{2}{\mathrm{Aut}(Q)} \end{array} \right). \nonumber
\end{equation}

\subsubsection{Rank 2 on $\mathbb{F}_{a}$}

In this section, we consider the more complicated case of rank 2 on $\mathbb{F}_{a}$ ($a \in \mathbb{Z}_{\geq 0}$). The fan of $\mathbb{F}_{a}$ is
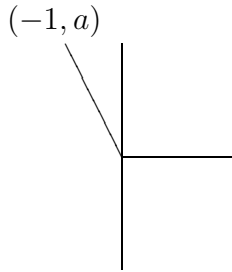
\begin{displaymath}
\xy
(0,0)*{} ; (0,15)*{} **\dir{-} ; (0,0)*{} ; (15,0)*{} **\dir{-} ; (0,0)*{} ; (0,-15)*{} **\dir{-} ; (0,0)*{} ; (-7.5,15)*{} **\dir{-} ; (-9,18)*{(-1,a)}, 
\endxy 
\end{displaymath}
so we obtain relations $D_{1} = D_{3}$ and $D_{4} = D_{2} + a D_{3}$ (Section 3.1). Defining $E := D_{1}$, $F := D_{2}$, the cohomology ring is given by
\begin{equation}
H^{2*}(\FF_a,\Z) \cong \mathbb{Z}[E,F]/(E^{2},F^{2}+aEF). \nonumber
\end{equation} 
A divisor $H = \alpha E + \beta F$ is ample if and only if $\beta > 0$ and $\alpha':=\alpha-a\beta>0$ \cite[Sect.~3.4]{Ful}. Fix such an ample divisor and an arbitrary first Chern class $c_{1} = f_{3}D_{3}+f_{4}D_{4} \in H^{2}(\FF_a,\mathbb{Z})$. By formula (\ref{ch. 2, eqn2}), the only interesting cases are $(f_{3},f_{4}) = (0,0),(1,0),(0,1),(1,1)$. 

\begin{corollary} \label{ch. 2, sect. 4, cor. 3}
Let $S = \mathbb{F}_{a}$, $H = \alpha D_{1} + \beta D_{2}$ an ample divisor, and $c_{1} = f_{3}D_{3} + f_{4}D_{4}$. Define $\lambda := \frac{\alpha}{\beta}$. The generating function $\prod_{k=1}^{\infty} (1-q^{k})^{8} \sum_{c_{2}} e( \M_{\FF_a}^{H}(2,c_{1},c_{2}) )q^{c_{2}}$ is given by
\begin{align}
- \sum_{(i,j,k,l) \in C_1} q^{\frac{1}{2}f_{3}f_{4}+\frac{a}{4}f_{4}^{2}+\frac{1}{2}j(i-\frac{a}{2}j)} + &2 \Big(  \sum_{(i,j,k,l) \in C_2} + \sum_{(i,j,k,l) \in C_3} \Big)q^{\frac{1}{2}f_{3}f_{4}+\frac{a}{4}f_{4}^{2}+\frac{1}{4}ij-\frac{1}{4}jk+\frac{1}{4}il+\frac{1}{4}kl-\frac{a}{4}l^{2}} \nonumber \\
+ &\Big( 2  \sum_{(i,j,k) \in C_4} +  \sum_{(i,j,k) \in C_5} + \sum_{(i,j,k) \in C_6} \Big) q^{\frac{1}{2}f_{3}f_{4}+\frac{a}{4}f_{4}^{2}+\frac{1}{2}j(i-\frac{a}{2}j)},  \nonumber
\end{align}
where $C_1, C_2, C_3 \subset \Z^4$, $C_5, C_6 \subset \Z^3$ are the following sets 
\begin{align*}
C_1 &:= \big\{(i,j,k,l) \in \Z^4 \ : \ 2 \ | \ f_{3}+i, \ 2 \ | \ f_{4}+j, \  2 \ | \ i+k, \ 2 \ | \ j+l, \ \lambda j = i, \ -j<l<j, \\ 
&\quad \quad -\lambda j + a(j+l) < k < \lambda j \big\}, \\
C_2 &:= \big\{(i,j,k,l) \in \Z^4 \ : \ 2 \ | \ f_{3} + i, \ 2 \ | \ f_{4} + j, \ 2 \ | \ i + k, \ 2 \ | \ j + l, \ k < \lambda l < i, \ l < j, \\
&\quad \quad -i - a(j - l) < k, \ -\lambda j < k \}, \\
C_3 &:= \big\{ (i,j,k,l) \in \Z^4 \ : \ 2 \ | \ f_{3} + i, \ 2 \ | \ f_{4} + j, \ 2 \ | \ i + k, \ 2 \ | \ j + l, \ k < \lambda l < i, \ l < j, \\
&\quad \quad -i + a(j + l) < k, \ -\lambda j + a(j + l) < k \big\}, \\
C_4 &:=  \big\{ (i,j,k) \in \Z^3 \ : \ 2 \ | \ f_{3} + i, \ 2 \ | \ f_{4} + j, \ 2 \ | \ j + k, \ i < \lambda j, \ \frac{a}{2}(j + k) < i, \\
&\quad \quad -\frac{i}{\lambda - a} + \frac{aj}{\lambda - a} < k < \lambda^{-1}i \big\}, \\
C_5 &:= \big\{ (i,j,k) \in \Z^3 \ : \ 2 \ | \ f_{3} + i, \ 2 \ | \ f_{4} + j, \ 2 \ | \ i + k, \ \lambda j < i, \ -\lambda j < k < \lambda j \big\}, \\
C_6 &:= \big\{ (i,j,k) \in \Z^3 \ : \ 2 \ | \ f_{3} + i, \ 2 \ | \ f_{4} + j, \ 2 \ | \ i + k, \ \lambda j < i, \ j > 0, \\
&\quad \quad -\lambda j + 2aj < k < \lambda j \big\}.
\end{align*}
\end{corollary}
\begin{proof}
Since $e(S)=4$ and $r=2$, the spaces $\D_{({\bf{v}},\bsdelta)}$ of Theorem \ref{ch. 2, sect. 3, thm. 1} are locally closed subsets of $(\PP^1)^4$. The only possibly non-empty quotients $\D_{({\bf{v}},\bsdelta)}^{s} / \mathrm{SL}(2,\C)$ occur for 
\begin{align*}
&\mathrm{all} \ v_i > 0 \ \mathrm{and \ all} \ \delta_i = 0, \\
&\mathrm{all} \ v_i > 0 \ \mathrm{and \ exactly \ one} \ \delta_i = 1, \\
&\mathrm{exactly  \ one} \ v_i = 0 \ \mathrm{and} \ \delta_j = 0 \ \mathrm{for \ all} \ j \neq i.
\end{align*}
The first line corresponds to moduli of four distinct points on $\PP^1$, \emph{or} moduli of four points on $\PP^1$ such that $p_1 = p_3$ and $p_1$, $p_2$, $p_4$ mutually distinct, \emph{or} moduli of four points on $\PP^1$ such that $p_2 = p_4$ and $p_1$, $p_2$, $p_3$ mutually distinct. This gives cases 1--3. The second line corresponds to moduli of four points on $\PP^1$ such that exactly two points coincide (the remaining possibilities: either $p_1 = p_2$, or $p_1 = p_4$, or $p_2 = p_3$, or $p_3 = p_4$). This gives cases 4--7. The third line corresponds to moduli of three distinct points on $\PP^1$. This gives cases 8--11. When non-empty, $e(\D_{({\bf{v}},\bsdelta)}^{s} / \mathrm{SL}(2,\C))$ is $-1$ in case one and $1$ in all other cases. Each of these eleven cases contributes one term to the generating function. Proceeding as in the previous section, we find that $\prod_{k=1}^{\infty} (1-q^{k})^{8} \sum_{c_{2}} e(\M_{\FF_a}^{H}(2,c_{1},c_{2})) q^{c_{2}}$ is equal to
\begin{align} \label{ch. 2, eqn3}
&-\sum_{{\footnotesize{ \begin{array}{c} v_{1}, v_{2}, v_{3}, v_{4} > 0 \ \mathrm{s.t.} \\ 2 \ | \ -f_{3}+v_{1}-a v_{2}+v_{3} \\ 2 \ | \ -f_{4}+v_{2}+v_{4} \\ \beta v_{1} < \alpha' v_{2}+\beta v_{3}+\alpha v_{4} \\ \alpha' v_{2}<\beta v_{1}+\beta v_{3}+\alpha v_{4} \\ \beta v_{3}<\beta v_{1}+\alpha' v_{2}+\alpha v_{4} \\ \alpha v_{4}<\beta v_{1}+\alpha' v_{2}+\beta v_{3} \end{array}   }}} q^{\frac{1}{2}f_{3}f_{4}+\frac{a}{4}f_{4}^{2}+\frac{1}{2}(v_{2}+v_{4})(v_{1}+\frac{a}{2}v_{2}+v_{3}-\frac{a}{2}v_{4})} \\
&+ \!\!\!\! \sum_{{\footnotesize{\begin{array}{c} v_{1}, v_{2}, v_{3}, v_{4} > 0 \ \mathrm{s.t.} \\ 2 \ | \ -f_{3}+v_{1}-a v_{2}+v_{3} \\ 2 \ | \ -f_{4}+v_{2}+v_{4} \\ \beta v_{1}+\alpha' v_{2}<\beta v_{3}+\alpha v_{4} \\ \beta v_{3}<\beta v_{1}+\alpha' v_{2}+\alpha v_{4} \\ \alpha v_{4}<\beta v_{1}+\alpha' v_{2}+\beta v_{3} \end{array}}}} \!\!\!\! q^{\frac{1}{2}f_{3}f_{4}+\frac{a}{4}f_{4}^{2}-\frac{1}{2}(v_{2}+v_{4})(v_{1}-\frac{a}{2} v_{2}+ v_{3}+\frac{a}{2} v_{4})+ v_{2} v_{3}+ v_{3} v_{4}+ v_{4} v_{1}} + \mathrm{5 \ similar \ terms} \nonumber \\
&+ \sum_{{\footnotesize{\begin{array}{c} v_{2}, v_{3}, v_{4} > 0 \ \mathrm{s.t.} \\ 2 \ | \ -f_{3}-a v_{2}+v_{3} \\ 2 \ | \ -f_{4}+v_{2}+v_{4} \\ \alpha' v_{2}<\beta v_{3}+\alpha v_{4} \\ \beta v_{3}<\alpha' v_{2}+\alpha v_{4} \\ \alpha v_{4}<\alpha' v_{2}+\beta v_{3} \end{array}}}} q^{\frac{1}{2}f_{3}f_{4}+\frac{a}{4}f_{4}^{2}+\frac{1}{2}(v_{2}+v_{4})(\frac{a}{2} v_{2}+v_{3}-\frac{a}{2} v_{4})} + \mathrm{3 \ similar \ terms}. \nonumber 
\end{align}
Next, we rewrite the first term and two of the next six terms of this expression. Specifically, we consider the term corresponding to all $p_i$'s mutually distinct and two of the terms corresponding to the cases where exactly two $p_i$'s coincide, namely the cases $p_1=p_3$ and $p_2=p_4$. For these three terms, we use the substitutions $i = v_{1}+v_{3}+a v_{2}$, $j = v_{2}+v_{4}$, $k = v_{1}-v_{3}+a v_{2}$ and $l = v_{2}-v_{4}$. After these substitutions, the terms combine to the first term of the corollary.

For the other four terms where exactly two $p_i$'s coincide, namely $p_1 = p_2$, $p_1=p_4$, $p_2=p_3$, $p_3=p_4$, we use the substitutions $i = v_{1}+v_{3}-av_{2}$, $j = v_{2}+v_{4}$, $k = v_{1}-v_{3}-a v_{2}$ and $l = - v_{2}+v_{4}$. This gives terms two and three of the corollary.

The last four terms of equation (\ref{ch. 2, eqn3}) can be rewritten as the last three terms of the corollary. For example, for the term corresponding to $v_1=0$, we use the substitutions $i = v_{3}+a v_{2}$, $j = v_{2}+v_{4}$ and $k = v_{2}-v_{4}$. The other three go similar.
\end{proof}

\begin{remark} \label{P^1xP^1}
Specializing the expression of Corollary \ref{ch. 2, sect. 4, cor. 3} to $a=0$ and setting $\lambda = \frac{\alpha}{\beta}$ gives
\begin{align}
&- \sum_{(i,j,k,l) \in C_{1}^{\prime}} q^{\frac{1}{2}f_{3}f_{4} + \frac{1}{2}ij} \nonumber + 4 \sum_{(i,j,k,l) \in C_{2}^{\prime}} q^{\frac{1}{2}f_{3}f_{4}+\frac{1}{4}ij-\frac{1}{4}jk+\frac{1}{4}il+\frac{1}{4}kl} + 2 \sum_{(i,j,k) \in C_{3}^{\prime} \cup C_{4}^{\prime}} q^{\frac{1}{2}f_{3}f_{4}+\frac{1}{2}ij},  \nonumber
\end{align}
where
\begin{align*}
C_{1}^{\prime} &:= \big\{ (i,j,k,l) \in \mathbb{Z}^4 \ : \  2 \ | \ f_{3}+i, \ 2 \ | \ f_{4}+j, \ 2 \ | \ i+k, \ 2 \ | \ j+l, \ \lambda j = i, \ -j<l<j, \\ 
&\quad \quad -\lambda j < k < \lambda j \big\}, \\
C_{2}^{\prime} &:= \big\{ (i,j,k,l) \in \mathbb{Z}^4 \ : \  2 \ | \ f_{3} + i, \ 2 \ | \ f_{4} + j, \ 2 \ | \ i + k, \ 2 \ | \ j + l, \ k < \lambda l < i, \ l < j, \\ 
&\quad \quad -i<k, \ -\lambda j < k \big\}, \\
C_{3}^{\prime} &:= \big\{ (i,j,k) \in \mathbb{Z}^3 \ : \ 2 \ | \ f_{3} + i, \ 2 \ | \ f_{4} + j, \ 2 \ | \ j + k, \ i < \lambda j, \ -\lambda^{-1} i < k < \lambda^{-1}i \big\}, \\
C_{4}^{\prime} &:= \big\{ (i,j,k) \in \mathbb{Z}^3 \ : \ 2 \ | \ f_{3} + i, \ 2 \ | \ f_{4} + j, \ 2 \ | \ i + k, \ \lambda j < i, \ -\lambda j < k < \lambda j \big\}. 
\end{align*}
Specializing to $\lambda=1$, i.e.~$H = D_1+D_2$, this expression can be simplified further. We consider the case $c_1 = D_3$, all other cases being similar. In this case, the generating function $\prod_{k=1}^{\infty} (1-q^{k})^{8} \sum_{c_{2}} e(\M_{\FF_a}^{H}(2,c_{1},c_{2})) q^{c_{2}}$ is given by
\begin{align}
&\sum_{m=1}^{\infty} \sum_{n=1}^{2m} \frac{4 q^{(2m+3)m-2mn+1}(q^{(2m+1)n}-q^{n^{2}})}{(1-q^{n})(q^{2m+1}-q^{n})} + \sum_{m=1}^{\infty} \frac{2(2m-1) q^{(2m-1)m}}{1-q^{2m-1}} + \sum_{m=1}^{\infty}\frac{4m q^{(2m+1)m}}{1-q^{2m}} \nonumber \\
&+ \sum_{m=1}^{\infty}\sum_{n=1}^{\infty}\sum_{p=1}^{2m-1} \frac{4 q^{(2m+1)m-2mp+1}((q^{n+p-1})^{p}-(q^{n+p-1})^{2m})}{q-q^{n+p}}. \qquad\qquad\qquad\qquad\qquad\qquad  \ \ \oslash \nonumber 
\end{align}
\end{remark}

\noindent \emph{Comparison to existing literature.} In \cite[Thm.~4.4]{Got2}, G\"ottsche gives an expression for generating functions of Hodge polynomials of moduli spaces of rank 2 $\mu$-stable torsion free sheaves on ruled surfaces $S$ with $-K_{S}$ effective. We consider this expression in the case $S = \mathbb{F}_{a}$. Among the toric divisors, $D_{1}$ is a fibre and $D_{2}$ is a section. Let $c_{1} = \epsilon D_{1} + D_{2}$ with $\epsilon \in \{0,1\}$, let $H$ be an ample divisor, and let $c_{2} \in H^4(\FF_a,\Z) \cong \mathbb{Z}$. Denote by $\M_{\FF_a}^{H,ss}(2,c_{1},c_{2})$ the moduli space of rank 2 Gieseker semistable torsion free sheaves on $\FF_a$ with Chern classes $c_{1}$, $c_{2}$. 
G\"ottsche and Qin \cite{Got2, Qin} have proved that the ample cone $C_{S} \subset \mathrm{Pic}(S) \otimes_{\mathbb{Z}} \mathbb{R}$ has a chamber/wall structure such that the moduli space $\M_{S}^{H,ss}(2,c_{1},c_{2})$ stays constant on each chamber. In our current example, the non-empty walls of type $(c_{1},c_{2})$ are  
\begin{equation*}
W^{\xi} = \{ x \in \mathrm{Pic}(\FF_a) \ \mathrm{ample} \ | \ x \cdot \xi = 0\},
\end{equation*}
where $\xi = (2n+\epsilon) D_{1}+(2m+1) D_{2}$ for any integers $m,n$ satisfying $m \geq 0$, $n <0$, $c_{2} - m(m+1)a + (2m+1)n + m \epsilon \geq 0$ \cite[Sect.~4]{Got2}. Elements $\frac{\alpha}{\beta} \in \mathbb{Q}_{>a}$ with $\alpha, \beta > 0$ coprime are in 1-1 correspondence with ample divisors $H = \alpha D_{1} + \beta D_{2}$ on $\FF_a$ with $\alpha$, $\beta$ coprime. Let $\Lambda$ be the collection of elements $\frac{\alpha}{\beta} \in \mathbb{Q}_{>a}$ with $\alpha$, $\beta$ coprime satisfying $\gcd(2, c_{1} \cdot (\alpha D_1 + \beta D_2)) = 1$. We refer to the complement $W = \mathbb{Q}_{>a} \setminus \Lambda$ as the collection of walls\footnote{The terminology ``wall'' might be slightly confusing in this context as $W$ lies dense in $\mathbb{Q}_{>a}$.}. The elements $\lambda \in \Lambda$ correspond to ample divisor $H$ for which there are no rank 2 strictly $\mu$-semistable torsion free sheaves with Chern class $c_{1}$ on $\FF_a$ \cite[Lem.~1.2.13, 1.2.14]{HL}. In this case $\M_{S}^{H,ss}(2,c_{1},c_{2}) = \M_{S}^{H}(2,c_{1},c_{2})$ for any $c_{2}$. The elements of $W$ are precisely the rational numbers corresponding to ample divisors lying on a wall of type $(c_{1},c_{2})$ for some $c_{2}$. For $H$ not on a wall as above \cite[Thm.~4.4]{Got2} gives
\begin{align}
\begin{split} \label{ch. 2, eqn4}  
&\sum_{c_{2}} e(\M_{\FF_a}^{H}(2,c_{1},c_{2})) q^{c_{2}} \\
&=\frac{1}{\prod_{k=1}^{\infty}(1-q^{k})^{8}} \sum_{(m, n) \in L(H)} \big( a + 2m a - 2(2 m + 2 n +\epsilon + 1) \big) q^{(m+1)m a-(2m+1)n - m \epsilon}, \\
&L(H) := \left\{ (m,n) \in \mathbb{Z}^{2} \ | \ m \geq 0, \ a - \lambda > \frac{2n + \epsilon}{2m + 1} \right\}. 
\end{split}
\end{align} 
Although G\"ottsche's formula (\ref{ch. 2, eqn4}) is equal to the formula of Corollary \ref{ch. 2, sect. 4, cor. 3}, it is not easy to obtain equality by direct manipulations. However, it is instructive to make expansions of both expressions for various values of $a$, $c_{1}$, $H$ (not on a wall). One finds a perfect agreement in such experiments.

\subsubsection{Wall-crossing for rank 2 on $\mathbb{F}_{a}$}

Theorem \ref{ch. 2, sect. 3, thm. 1} also allows one to study the dependence on choice of polarization. This leads to wall-crossing formulae. We illustrate this in the case of rank 2 sheaves on $\FF_a$. We start with a few definitions. Denote by $\mathbb{Z}(\!(q)\!)$ the ring of formal Laurent series over $\Z$. For all values $\lambda \in \mathbb{Q}_{>a}$ of the stability parameter, the expression of Corollary \ref{ch. 2, sect. 4, cor. 3} is a formal Laurent series. Therefore, we can see the expression of Corollary \ref{ch. 2, sect. 4, cor. 3} as a map $\mathbb{Q}_{>a} \longrightarrow \mathbb{Z}(\!(q)\!)$. We define the following notion of limit.
\begin{definition} \label{ch. 2, sect. 4, def. 1}
Let $a \in \mathbb{Z}_{\geq 0}$ and let $F : \mathbb{Q}_{>a} \longrightarrow \mathbb{Z}(\!(q)\!)$, $\lambda \mapsto F(\lambda)$ be a map. Let $\lambda_{0} \in \mathbb{Q}_{>a}$ and let $F_{0} \in \mathbb{Z}(\!(q)\!)$. We define
\begin{equation*}
\lim_{\epsilon, \epsilon' \searrow 0 } \left( F(\lambda_{0} + \epsilon) - F(\lambda_{0} - \epsilon') \right) = F_{0}
\end{equation*}
whenever for any $N \in \mathbb{Z}$ there exist $\epsilon, \epsilon' \in \mathbb{Q}_{>0}$ such that $a < \lambda_{0} - \epsilon'$ and
\begin{equation*}
F(\lambda_{0} + \epsilon) - F(\lambda_{0} - \epsilon') = F_{0} \mod q^N. 
\end{equation*} 
Note that if the limit exists, then it is unique. We refer to the expression
\begin{equation*}
\lim_{\epsilon, \epsilon' \searrow 0 } \left( F(\lambda_{0} + \epsilon) - F(\lambda_{0} - \epsilon') \right) = F_{0}
\end{equation*}
as an \emph{infinitesimal wall-crossing formula}. \hfill $\oslash$
\end{definition}

\noindent Applying this notion of limit to the expression of Remark \ref{P^1xP^1} gives the following result. 
\begin{corollary} \label{ch. 2, sect. 4, cor. 6}
Let $S= \mathbb{P}^{1} \times \mathbb{P}^{1}$, let $H = \alpha_{0} D_{1} + \beta_{0} D_{2}$ be an ample divisor, and suppose without loss of generality that $\gcd(\alpha_{0},\beta_{0})=1$. Let $c_{1} = f_{3}D_{3} + f_{4}D_{4} \in H^{2}(S,\mathbb{Z})$. Defining $\lambda_{0} = \frac{\alpha_{0}}{\beta_{0}}$, we have
\begin{align*}
&\lim_{\epsilon, \epsilon' \searrow 0} \prod_{k=1}^{\infty} (1-q^{k})^{8} \left( \sum_{c_{2}} e( \M_{\mathbb{P}^{1} \times \mathbb{P}^{1}}^{\lambda_{0}+\epsilon}(2,c_{1},c_{2}) )q^{c_{2}} - \sum_{c_{2}} e( \M_{\mathbb{P}^{1} \times \mathbb{P}^{1}}^{\lambda_{0}-\epsilon'}(2,c_{1},c_{2}) )q^{c_{2}} \right) \\ 
&= 4 \Big( \sum_{(i,j,k) \in C_{1}^{\prime \prime}} - \sum_{(i,j,k) \in C_{2}^{\prime \prime}} \Big) q^{\frac{1}{2}f_{3}f_{4}+\frac{1}{4}ij - \frac{\lambda_{0}}{4}jk + \frac{1}{4}ik + \frac{\lambda_{0}}{4}k^{2}} \\
&\quad + 4 \Big( \sum_{(i,j,k) \in C_{3}^{\prime \prime}} - \sum_{(i,j,k) \in C_{4}^{\prime \prime}} \Big) q^{\frac{1}{2}f_{3}f_{4}+\frac{1}{4}ij - \frac{\lambda_{0}}{4}jk + \frac{1}{4}ik + \frac{\lambda_{0}}{4}k^{2}} \displaybreak \\
&\quad+ \sum_{(i,j) \in C_{5}^{\prime \prime}} 2 q^{\frac{1}{2}f_{3}f_{4}+\frac{\lambda_{0}}{2}i^{2}} -  \sum_{(i,j) \in C_{6}^{\prime \prime}} 2 q^{\frac{1}{2}f_{3}f_{4}+\frac{\lambda_{0}^{-1}}{2}i^{2}} \\
&\quad+ \sum_{(i,j) \in C_{7}^{\prime \prime}} 4 q^{\frac{1}{2}f_{3}f_{4}+\frac{\lambda_{0}^{-1}}{2}ij} - \sum_{(i,j) \in C_{8}^{\prime \prime}} 4 q^{\frac{1}{2}f_{3}f_{4}+\frac{\lambda_{0}}{2}ij}, 
\end{align*}
where
\begin{align*}
&C_{1}^{\prime \prime} := \big\{ (i,j,k) \in \mathbb{Z}^3 \ : \ \beta_{0} \ | \ k, \ 2 \ | \ f_{3} + i, \ 2 \ | \ f_{4} + j, \ 2 \ | \ i + \lambda_{0}k, \ 2 \ | \ j+k,\\ 
&\qquad \qquad 0<\lambda_{0}k<i, \ 0<k<j \big\}, \\
&C_{2}^{\prime \prime} := \big\{ (i,j,k) \in \mathbb{Z}^3 \ : \  \beta_{0} \ | \ k, \ 2 \ | \ f_{3} + i, \ 2 \ | \ f_{4} + j, \ 2 \ | \ i + \lambda_{0}k, \ 2 \ | \ j+k, \\ 
&\qquad \qquad -i<\lambda_{0}k<0, \ -j<k<0 \big\}, \\
&C_{3}^{\prime \prime} := \big\{ (i,j,k) \in \mathbb{Z}^3 \ : \  \beta_{0} \ | \ k, \ 2 \ | \ f_{3} + i, \ 2 \ | \ f_{4} + k, \ 2 \ | \ i + \lambda_{0}k, \ 2 \ | \ j+k, \\ 
&\qquad \qquad -k<j<k, \ \lambda_{0}k<i \big\}, \\
&C_{4}^{\prime \prime} := \big\{ (i,j,k) \in \mathbb{Z}^3 \ : \ \beta_{0} \ | \ k, \ 2 \ | \ f_{3} + \lambda_{0}k, \ 2 \ | \ f_{4} + j, \ 2 \ | \ i + \lambda_{0}k, \ 2 \ | \ j+k, \\ 
&\qquad \qquad -\lambda_{0}k<i<\lambda_{0}k, \ k<-j \big\}, \\
&C_{5}^{\prime \prime} := \big\{(i,j) \in \mathbb{Z}^2 \ : \ \beta_{0} \ | \ i, \ 2 \ | \ f_{3}+\lambda_{0}i, \ 2 \ | \ f_{4}+i, \ 2 \ | \ i+j, \ -i<j<i \big\}, \\
&C_{6}^{\prime \prime} := \big\{ (i,j) \in \mathbb{Z}^2 \ : \ \alpha_{0} \ | \ i, \ 2 \ | \ f_{4}+\lambda_{0}^{-1}i, \ 2 \ | \ f_{3}+i, \ 2 \ | \ i+j, \ -i<j<i \big\}, \\
&C_{7}^{\prime \prime} := \big\{ (i,j) \in \mathbb{Z}^2 \ : \ \alpha_{0} \ | \ j, \ 2 \ | \ f_{4}+\lambda_{0}^{-1}j, \ 2 \ | \ f_{3}+i, \ 2 \ | \ i+j, \ 0<j<i \big\}, \\
&C_{8}^{\prime \prime} := \big\{ (i,j) \in \mathbb{Z}^2 \ : \ \beta_{0} \ | \ j, \ 2 \ | \ f_{3}+\lambda_{0}j, \ 2 \ | \ f_{4}+i, \ 2 \ | \ i+j, \ 0<j<i \big\}.
\end{align*}
\end{corollary} 

Roughly speaking, the formula of the previous corollary is obtained from all possible ways of changing an inequality involving $\lambda$ in the formula of Remark \ref{P^1xP^1} into an equality and summing these terms with appropriate signs. The expression of the previous corollary can only be non-zero when $2 \ | \ \alpha_{0} f_{4} + \beta_{0} f_{3}$ or, equivalently, $H$ lies on a wall. \\

\noindent \emph{Comparison to existing literature.} It is easy to derive a nice infinitesimal wall-crossing formula from G\"ottsche's formula (\ref{ch. 2, eqn4}). Let $c_{1} = \epsilon D_{1} + D_{2}$ ($\epsilon \in \{0,1\}$) and $\lambda_{0} = \frac{\alpha_{0}}{\beta_{0}} \in \Q_{>a}$ arbitrary (i.e.~corresponding to any ample divisor $H = \alpha_0 D_1 + \beta_0 D_2$ with $\beta_0 > 0$, $\alpha_0 > a \beta_0$, and $\gcd(\alpha_0,\beta_0) = 1$). Using Definition \ref{ch. 2, sect. 4, def. 1}, one obtains
\begin{align} 
\begin{split} \label{ch. 2, eqn5}
&\lim_{\epsilon, \epsilon' \searrow 0} \prod_{k=1}^{\infty}(1-q^{k})^{8} \left( \sum_{c_{2}} e( \M_{\mathbb{F}_{a}}^{\lambda_{0}+\epsilon}(2,c_{1},c_{2}))q^{c_{2}} - \sum_{c_{2}} e( \M_{\mathbb{F}_{a}}^{\lambda_{0}-\epsilon'}(2,c_{1},c_{2}))q^{c_{2}} \right) \\
&= \sum_{{\footnotesize{\begin{array}{c} m \in \mathbb{Z}_{\geq 1} \ \mathrm{s.t.} \\ \frac{1}{2} (\lambda_{0}-a)(2m-1) - \frac{1}{2} \epsilon \in \mathbb{Z} \end{array}}}}  2\left(1+\frac{a}{2}-\lambda_{0}\right)\left(2m-1\right)q^{\frac{1}{2}(\lambda_{0}-\frac{a}{2})(2m-1)^{2}-\frac{1}{4}a+\frac{1}{2}\epsilon}. 
\end{split}
\end{align}
Since the complement of all walls $\Lambda \subset \Q_{>a}$ lies dense, strictly $\mu$-semistables do not play a role in this formula.

We can also derive equation (\ref{ch. 2, eqn5}) using Joyce's machinery for wall-crossing of motivic invariants counting (semi)stable objects in an abelian category \cite{Joy2}. Joyce gives a wall-crossing formula for virtual Poincar\'e polynomials of moduli stacks of Gieseker semistable torsion free sheaves on an arbitrary nonsingular complete surface $S$ with $-K_{S}$ nef \cite[Thm.~6.21]{Joy2}\footnote{The cited theorem also holds for $\mu$-stability instead of Gieseker stability.}. For $S = \FF_a$, these are $\mathbb{P}^{1} \times \mathbb{P}^{1}$, $\mathbb{F}_{1}$, $\mathbb{F}_{2}$. Nevertheless, we apply the formula of \cite[Thm.~6.21]{Joy2} to any $S = \FF_a$ keeping $a$ arbitrary. Let $c_{1} = f_{3}D_{3}+f_{4}D_{4} \in H^{2}(\FF_a,\mathbb{Z})$ and $\lambda_0 = \frac{\alpha_0}{\beta_0}$ as before. Part of Joyce's philosophy is to study wall-crossing phenomena for motivic invariants of moduli \emph{stacks} instead of moduli \emph{schemes} (coming from GIT as in \cite[Ch.~4]{HL}). Keeping track of the stabilizers gives nice wall-crossing formulae. In this paper we are interested in Euler characteristics of moduli schemes (coming from GIT as in \cite[Ch.~4]{HL}), so we first make a connection between the two. 

For any smooth complete surface $S$, polarization $H$, $r > 0$, and Chern classes $c_1$, $c_2$ let $\M_{S}^{H}(r,c_{1},c_{2})$ be the coarse moduli scheme of rank $r$ $\mu$-stable torsion free sheaves on $S$ with Chern classes $c_1$, $c_2$ as before. Let $\mathrm{Obj}_{s}^{\mathrm{ch}}(\mu)$ be the Artin stack of $\mu$-stable torsion free sheaves on $S$ with total Chern character $\mathrm{ch} = (r,c_1,\frac{1}{2}(c_{1}^{2}-2c_2))$ \cite{Joy2}. Denote the virtual Poincar\'e polynomial by $P(\cdot,z)$. Joyce proves one can uniquely extend the definition of virtual Poincar\'e polynomial to Artin stacks of finite type over $\mathbb{C}$ with affine geometric stabilizers if one requires
\[
P([Y/G],z) = P(Y,z) / P(G,z)
\]
for any \emph{special} algebraic group $G$ acting regularly on a quasi-projective variety $Y$ \cite[Thm.~4.10]{Joy1}. We claim
\begin{equation} \label{ch. 2, eqn6}
e(\M_{S}^{H}(r,c_{1},c_{2})) = \lim_{z \rightarrow -1} \left( (z^{2}-1) P(\mathrm{Obj}_{s}^{\mathrm{ch}}(\mu), z) \right).
\end{equation} 
This equation can be proved as follows. Recall that $\M_{S}^{H}(r,c_{1},c_{2})$ is constructed as a geometric quotient $\pi : R^{s} \longrightarrow \M_{S}^{H}(r,c_{1},c_{2})$, where $R^s$ is an open subset of some Quot scheme with an action of $\mathrm{PGL}(n,\mathbb{C})$ for some $n$ \cite[Ch.~4]{HL}. In fact, $\pi$ is a principal $\mathrm{PGL}(n,\mathbb{C})$-bundle \cite[Cor.~4.3.5]{HL} and we have isomorphisms of stacks \cite[Prop.~3.3]{Gom}
\begin{equation*}
\M_{S}^{H}(r,c_{1},c_{2}) \cong [ R^{s} / \mathrm{PGL}(n,\mathbb{C}) ], \ \mathrm{Obj}_{s}^{\mathrm{ch}}(\mu) \cong [ R^{s} / \mathrm{GL}(n,\mathbb{C}) ].
\end{equation*}
The difficulty is that $\mathrm{PGL}(n,\mathbb{C})$ is in general \emph{not} special. Let $(\C^*)^n \leq \mathrm{GL}(n,\C)$ be the subgroup of diagonal matrices. Define $\mathrm{P}(\mathbb{C}^{*})^{n} = (\mathbb{C}^{*})^{n} / \mathbb{C}^{*} \cdot \mathrm{id}$, where $\mathrm{id}$ is the $n \times n$ identity matrix, and consider the geometric quotient $R^{s} / \mathrm{P}(\mathbb{C}^{*})^{n}$. We obtain a morphism 
\[
R^{s} / \mathrm{P}(\mathbb{C}^{*})^{n} \longrightarrow R^{s} / \mathrm{PGL}(n,\mathbb{C}),
\]
and all fibres over closed points are isomorphic to $F = \mathrm{PGL}(n, \mathbb{C}) / \mathrm{P}(\mathbb{C}^{*})^{n}$. We obtain 
\begin{align*}
e\left( \M_{S}^{H}(r,c_{1},c_{2}) \right) &= \frac{e\left( R^{s} / \mathrm{P}(\mathbb{C}^{*})^{n} \right)}{e(F)} = \frac{e\left( R^{s} / \mathrm{P}(\mathbb{C}^{*})^{n} \right)}{n!} = \lim_{z \rightarrow -1 }  \frac{P(R^{s},z)}{n!(z^{2}-1)^{n-1}} \\
&= \lim_{z \rightarrow -1 } \frac{(z^{2}-1)P(R^{s},z)}{P(\mathrm{GL}(n,\mathbb{C}),z)} \cdot \frac{(z^{2})^{\frac{n(n-1)}{2}}\prod_{k=1}^{n}((z^{2})^{k}-1)}{n!(z^{2}-1)^{n}}, 
\end{align*}
where we use \cite[Thm.~2.4]{Joy2} and \cite[Lem.~4.6]{Joy1}. Using
\begin{align*}
&\mathrm{lim}_{z \rightarrow -1} \frac{(z^{2})^{\frac{n(n-1)}{2}}\prod_{k=1}^{n}((z^{2})^{k}-1)}{(z^{2}-1)^{n}} = n!, \\
&\frac{P(R^s,z)}{P(\mathrm{GL}(n,\C),z)} = P([R^{s}/\mathrm{GL}(n,\mathbb{C})],z),
\end{align*}
we obtain formula (\ref{ch. 2, eqn6}). 

Back to $S = \FF_a$, using equation (\ref{ch. 2, eqn6}) and \cite[Thm.~6.21]{Joy2} a somewhat lengthy computation gives 
\begin{align} 
&\lim_{\epsilon, \epsilon' \searrow 0} \prod_{k=1}^{\infty}(1-q^{k})^{8} \left( \sum_{c_{2}} e( \M_{\FF_a}^{\lambda_{0}+\epsilon}(2,c_{1},c_{2}))q^{c_{2}} - \sum_{c_{2}} e( \M_{\FF_a}^{\lambda_{0}-\epsilon'}(2,c_{1},c_{2}))q^{c_{2}} \right) \label{ch. 2, eqn7} \\
&= \!\!\!\!\!\!\!\!\!\! \sum_{{\footnotesize{\begin{array}{c} m \in \mathbb{Z}_{>\frac{1}{2}f_{4}} \ \mathrm{s.t.} \\ \frac{1}{2} (\lambda_{0}-a)(2m-f_{4}) - \frac{1}{2}(f_{3}+af_{4}) \in \mathbb{Z} \end{array}}}} \!\!\!\!\!\!\!\!\!\! \!\!\!\!\!\!\!\!\!\! 2\left(1+\frac{a}{2}-\lambda_{0}\right)\left(2m-f_{4}\right)q^{\frac{1}{2}(\lambda_{0}-\frac{a}{2})(2m-f_{4})^{2}-\frac{1}{4}af_{4}^{2}+\frac{1}{2}(f_{3}+af_{4})f_{4}}. \nonumber
\end{align}
Note that \cite[Thm.~6.21]{Joy2} is a wall-crossing formula for Artin stacks of \emph{semi}stable objects, whereas we have been dealing with Artin stacks of stable objects only. In the cases $f_{3} \neq 0 \mod 2$ or $f_{4} \neq 0 \mod 2$, the complement of all walls, i.e.~$\Lambda \subset \Q_{>a}$, lies dense, so strictly $\mu$-semistables do not play a role in the above formula\footnote{In the case $f_3=f_4=0 \mod 2$, we have $\Lambda = \varnothing$. However, for $r=2$ and fixed $c_1, c_2$ one can show that $\mathrm{Obj}_{ss}^{\mathrm{ch}}(\mu) \setminus \mathrm{Obj}_{s}^{\mathrm{ch}}(\mu)$ is the same for any polarization not on a wall of type $(c_1,c_2)$. Therefore, formula (\ref{ch. 2, eqn7}) also holds in this case, because strictly $\mu$-semistables on either side of a wall cancel (compare \cite[Thm.~2.9]{Got2}).}. Note that equations (\ref{ch. 2, eqn5}) and (\ref{ch. 2, eqn7}) are consistent. In fact, they are even consistent in the case $a > 2$ suggesting \cite[Thm.~6.21]{Joy2} holds more generally. 

We now proved
expressions (\ref{ch. 2, eqn5}) and (\ref{ch. 2, eqn7}) obtained from G\"ottsche's and Joyce's work are equal to the wall-crossing formulae obtained from Corollary \ref{ch. 2, sect. 4, cor. 3} (e.g.~Corollary \ref{ch. 2, sect. 4, cor. 6} when $a=0$). This is by no means clear from direct manipulations of the expressions. It is instructive to make expansions to a certain order for various values of $a$, $\lambda_0$, $f_3$, $f_4$ and verify consistency. Similar to Remark \ref{P^1xP^1}, the wall-crossing formula of Corollary \ref{ch. 2, sect. 4, cor. 6} can be simplified for specific values of $\lambda_0$. We will not write down the explicit expressions.

\subsection{Rank 3 on $\mathbb{P}^{2}$}

We now apply\footnote{During the final preparations of the first version of this paper, the author found out about recent independent work of Weist \cite{Wei}, which also computes the case $r=3$ and $S = \mathbb{P}^{2}$ using techniques of toric geometry and quivers. Weist has communicated to the author that his results are consistent with the expansions given at the end of this section.} Theorem \ref{ch. 2, sect. 3, thm. 1} to the case $r=3$ and $S = \PP^2$. Similar computations can be done in the case $r=3$ and $S = \mathbb{F}_{a}$, but the formulae become (even) lengthier.  

Let $c_1 = f H$, where $H$ is the hyperplane class. Consider the expression of Theorem \ref{ch. 2, sect. 3, thm. 1}. Let $v_i := v_{i,1}$, $w_i := v_{i,2}$, $p_i := p_{i,1}$, and $q_i := p_{i,2}$. Moreover, let $v:= \sum_i v_i$ and $w:= \sum_i w_i$. For $v_i$, $w_i$ all positive and any choice of $\bsdelta$, we have
\[
\D_{({\bf{v}}, \bsdelta)}^{s} \subset \{(p_1, p_2, p_3, q_1, q_2, q_3) \ : \ p_i \subset q_i \ \forall i\} \subset \mathrm{Gr}(1,3)^3 \times \mathrm{Gr}(2,3)^3 \cong (\mathbb{P}^{2})^{3} \times (\mathbb{P}^{2 *})^{3}.
\]
Suppose all $\delta_{a,i} = 0$. Then $\D_{({\bf{v}}, \bsdelta)}^{s}$ is empty unless
\[
({\bf{v}},{\bf{w}}) := (v_1,v_2,v_3,w_1,w_2,w_3) \in C_1 \cup C_2,
\]
where
\begin{align*}
&C_1 := \big\{ ({\bf{v}},{\bf{w}}) \in \Z_{>0}^{6} \ : \ 3 \ | \ -f + v + 2w, \\ 
&\qquad \qquad v_{i} + 2w_{i} < 2v_{j} + 2v_{k} + w_{j} + w_{k}, \ w_{i} + 2v_{i} < 2w_{j} + 2w_{k} + v_{j} + v_{k}, \\ 
&\qquad \qquad v_{i} + v_{j} < 2v_{k} + w, \ w_{i} + w_{j} < 2w_{k} + v \ \forall \{i,j,k\} = \{1,2,3\} \big\}, \displaybreak \\
&C_2 := \big\{ ({\bf{v}},{\bf{w}}) \in \Z_{>0}^{6} \ : \ 3 \ | \ -f + v + 2w, \\ 
&\qquad \qquad v_{i} + 2w_{i} < 2v_{j} + 2v_{k} + w_{j} + w_{k}, \ w_{i} + 2v_{i} < 2w_{j} + 2w_{k} + v_{j} + v_{k}, \\ 
&\qquad \qquad v < w , \ w_{i} + w_{j} < 2w_{k} + v \ \forall \{i,j,k\} = \{1,2,3\} \big\}.
\end{align*}
The notation ``for all $\{i,j,k\}=\{1,2,3\}$'' means ``for all $i \in \{1,2,3\}$, $j \in \{1,2,3\} \setminus \{i\}$, and $k \in \{1,2,3\} \setminus \{i,j\}$''. For $({\bf{v}},{\bf{w}}) \in C_1$, $\D_{({\bf{v}}, \bsdelta)}^{s}$ is equal to the configuration space of $(p_1, p_2, p_3, q_1, q_2, q_3)$, where $q_i \subset \PP^2$ are lines such that $q_{1} \cap q_{2}$, $q_{2} \cap q_{3}$, $q_{3} \cap q_{1}$ are mutually distinct points, $p_i \subset q_i$ are points not equal to $q_{1} \cap q_{2}$, $q_{2} \cap q_{3}$, $q_3 \cap q_1$ and are not colinear. We denote this space pictorially by 
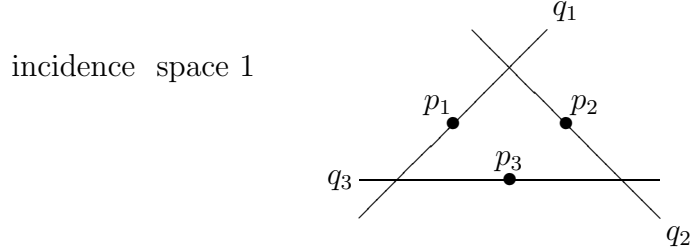
\begin{displaymath}
\xy
(-30,20)*{\textrm{incidence \ space} \ 1} ; (0,0)*{} ; (25,25)*{} **\dir{-} ; (15,25)*{} ; (40,0)*{} **\dir{-} ; (0,5)*{} ; (40,5)*{} **\dir{-} ; (12.5,12.5)*{\bullet} ; (27.5,12.5)*{\bullet} ; (20,5)*{\bullet} ; (10.5,15)*{p_{1}} ; (30,15)*{p_{2}} ; (20,7.5)*{p_{3}} ; (27.5,27.5)*{q_{1}} ; (42.5,-2.5)*{q_{2}} ; (-2.5,5)*{q_{3}}
\endxy 
\end{displaymath}
After taking the quotient by $\mathrm{SL}(3,\C)$, one obtains a space with Euler characteristic $-1$. This can be seen by using that for any four points $x_1$, $x_2$, $x_3$, $x_4$ of $\PP^2$, no three of which are colinear, there exists an element $g \in \mathrm{SL}(3,\C)$ mapping them to $(1:0:0)$, $(0:1:0)$, $(0:0:1)$, and $(1:1:1)$ respectively. Moreover, $g$ is unique up to multiplication by a 3rd root of unity. For $({\bf{v}},{\bf{w}}) \in C_2$, the incidence space is
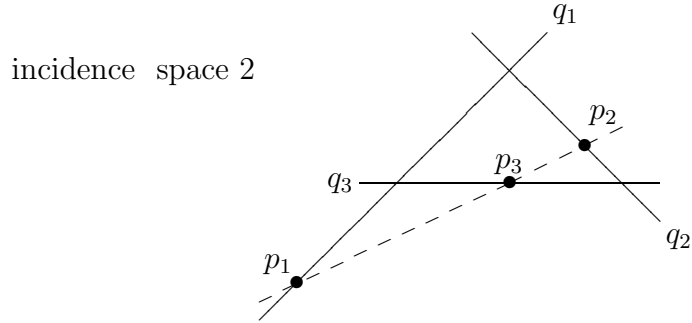
\begin{displaymath}
\xy
(-30,20)*{\textrm{incidence \ space} \ 2} ; (0,0)*{} ; (25,25)*{} **\dir{-} ; (15,25)*{} ; (40,0)*{} **\dir{-} ; (0,5)*{} ; (40,5)*{} **\dir{-} ; (12.5,12.5)*{} ; (27.5,12.5) ; (20,5) ; (10.5,15)*{} ; (27.5,27.5)*{q_{1}} ; (42.5,-2.5)*{q_{2}} ; (-2.5,5)*{q_{3}} ; (35,12.5)*{} ; (20,5)*{} **\dir{--} ; (20,5)*{} ; (-13.25,-10.75)*{} **\dir{--} ; (0,0)*{} ; (-13.25,-13.25) **\dir{-} ; (30,10)*{\bullet} ; (20,5)*{\bullet} ; (-8.25,-8.25)*{\bullet}; ; (20,7.5)*{p_{3}} ; (32.5,14)*{p_{2}} ; (-10.75,-5.75)*{p_{1}}
\endxy 
\end{displaymath}
where the dashed lines means $p_{1}$, $p_{2}$, $p_{3}$ are colinear. After taking the quotient by $\mathrm{SL}(3,\C)$, one obtains a reduced point. The contribution of these two incidence spaces to the generating function $q^{-\frac{1}{2} f^{2}} \prod_{k=1}^{\infty} (1-q^{k})^{9} \sum_{c_{2}} e(\M_{\PP^2}(3,c_{1},c_{2}))q^{c_{2}}$ is
\begin{align*}
&\Big(- \sum_{({\bf{v}},{\bf{w}}) \in C_1} + \sum_{({\bf{v}},{\bf{w}}) \in C_2} \Big)q^{Q_1({\bf{v}},{\bf{w}})}, \ \mathrm{where} \\
&Q_1({\bf{v}},{\bf{w}}) := - \frac{1}{18}(-f-2v-w)^{2} - \frac{1}{18}(-f+v-w)^{2} - \frac{1}{18}(-f+v+2w)^{2} \\
&\qquad \qquad \quad \ + \sum_{i<j}(v_i + w_i)(v_j+w_j).
\end{align*}

Similarly, other choices of $({\bf{v}},{\bf{w}},\bsdelta)$ give rise to other systems of inequalities and corresponding incidence spaces. We list all other incidence spaces which contribute.

\begin{displaymath}
\xy
(-40,20)*{\textrm{incidence \ space} \ 3} ; (0,0)*{} ; (25,25)*{} **\dir{-} ; (0,25)*{} ; (25,0)*{} **\dir{-} ; (-5,12.5)*{} ; (30,12.5)*{} **\dir{-} ; (27.5,-2.5)*{q_{1}} ; (-2.5,-2.5)*{q_{2}} ; (32.5,12.5)*{q_{3}} ; (6.25, 18.75)*{\bullet} ; (8.75, 21.25)*{p_{1}} ; (18.75, 18.75)*{\bullet} ; (16.25, 21.25)*{p_{2}} ; (22.5,12.5)*{\bullet} ; (22.5,15)*{p_{3}} 
\endxy
\end{displaymath}
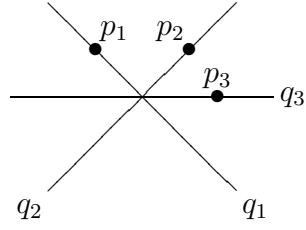
\begin{displaymath}
\xy
(-35,20)*{\textrm{incidence \ spaces} \ 4\textrm{--}9} ; (0,0)*{} ; (25,25)*{} **\dir{-} ; (15,25)*{} ; (40,0)*{} **\dir{-} ; (0,5)*{} ; (40,5)*{} **\dir{-} ; (5,5)*{\bullet} ; (27.5,12.5)*{\bullet} ; (20,5)*{\bullet} ; (2.5,7.5)*{p_{i}} ; (30,15)*{p_{j}} ; (20,7.5)*{p_{k}} ; (27.5,27.5)*{q_{i}} ; (42.5,-2.5)*{q_{j}} ; (-2.5,5)*{q_{k}}
\endxy 
\end{displaymath}
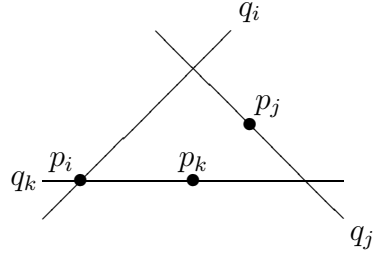
\begin{displaymath}
\xy
(-30,20)*{\textrm{incidence \ spaces} \ 10, 11, 12} ; (0,0)*{} ; (25,25)*{} **\dir{-} ; (15,25)*{} ; (40,0)*{} **\dir{-} ; (0,5)*{} ; (40,5)*{} **\dir{-} ; (27.5,12.5)*{\bullet} ; (20,5)*{\bullet} ; (30,15)*{p_{j}} ; (20,7.5)*{p_{k}} ; (27.5,27.5)*{q_{i}} ; (42.5,-2.5)*{q_{j}} ; (-2.5,5)*{q_{k}}
\endxy 
\end{displaymath}
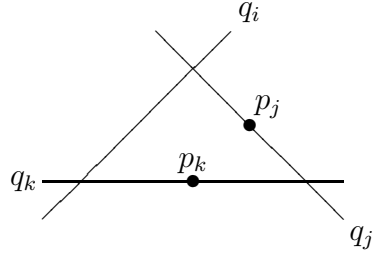
\begin{displaymath}
\xy
(-30,20)*{\textrm{incidence \ spaces} \ 13, 14, 15} ; (0,0)*{} ; (25,25)*{} **\dir{-} ; (15,25)*{} ; (40,0)*{} **\dir{-} ; (12.5,12.5)*{\bullet} ; (27.5,12.5)*{\bullet} ; (20,5)*{\bullet} ; (10.5,15)*{p_{i}} ; (30,15)*{p_{j}} ; (20,7.5)*{p_{k}} ; (27.5,27.5)*{q_{i}} ; (42.5,-2.5)*{q_{j}} ; (-2.5,5)*{}
\endxy 
\end{displaymath}
for all $\{i,j,k\} = \{1,2,3\}$. For incidence spaces 3 and 13--15, $p_{1}$, $p_{2}$, $p_{3}$ are not colinear. The incidence spaces 4--9 all give the same contribution to the generating function. This also holds for incidence spaces 10, 11, 12 as well as incidence spaces 13, 14, 15. The final answer is
\begin{align*}
&q^{-\frac{1}{2} f^{2}} \prod_{k=1}^{\infty} (1-q^{k})^{9} \sum_{c_{2}} e(\M_{\PP^2}(3,c_{1},c_{2}))q^{c_{2}} \\
&=\Big(- \sum_{({\bf{v}},{\bf{w}}) \in C_1} + \sum_{({\bf{v}},{\bf{w}}) \in C_2} + \sum_{({\bf{v}},{\bf{w}}) \in C_3} \Big)q^{Q_1({\bf{v}},{\bf{w}})} \\
&\quad \ + \sum_{({\bf{v}},{\bf{w}}) \in C_4} 6q^{Q_2({\bf{v}},{\bf{w}})} + \sum_{({\bf{v}},{\bf{w}}) \in C_5} 3q^{Q_1({\bf{v}},{\bf{w}})} + \sum_{({\bf{v}},{\bf{w}}) \in C_6} 3q^{Q_1({\bf{v}},{\bf{w}})}, 
\end{align*}
where
\begin{align*}
&Q_1({\bf{v}},{\bf{w}}) \ \mathrm{defined \ above}, \\ 
&Q_2({\bf{v}},{\bf{w}}) := Q_1({\bf{v}},{\bf{w}}) -v_1 w_3,
\end{align*}
\begin{align*}
&C_1, C_2 \ \mathrm{defined \ above}, \\
&C_3 := \big\{ ({\bf{v}},{\bf{w}}) \in \Z_{>0}^{6} \ : \ 3 \ | \ -f + v + 2w, \\ 
&\qquad \qquad v_{i} + 2w_{i} < 2v_{j} + 2v_{k} + w_{j} + w_{k}, \ w_{i} + 2v_{i} < 2w_{j} + 2w_{k} + v_{j} + v_{k}, \\ 
&\qquad \qquad w < v , \ v_{i} + v_{j} < 2v_{k} + w \ \forall \{i,j,k\} = \{1,2,3\} \big\}, \\
&C_4 := \big\{  ({\bf{v}},{\bf{w}}) \in \Z_{>0}^{6} \ : \ 3 \ | \ -f + v + 2w, \ v_{1} + 2w_{1} < 2v_{2} + 2v_{3} + w_{2} + w_{3}, \\ 
&\qquad \qquad v_{2} + 2w_{2} < 2v_{1} + 2v_{3} + w_{1} + w_{3}, \ w_{2} + 2v_{2} < 2w_{1} + 2w_{3} + v_{1} + v_{3}, \\ 
&\qquad \qquad w_{3} + 2v_{3} < 2w_{1} + 2w_{2} + v_{1} + v_{2}, \  v_{1} + v_{2} < 2v_{3} + w, \ v_{2} + v_{3} < 2v_{1} + w,  \\ 
&\qquad \qquad w_{1} + w_{2} < 2w_{3} + v, \ w_{2} + w_{3} < 2w_{1} + v, \ v_{1} + v_{3} + 2w_{3} < 2v_{2} + w_{1} + w_{2} \\ 
&\qquad \qquad w_{1} + w_{3} + 2v_{1} < 2w_{2} + v_{2} + v_{3} \big\}, \\
&C_5 := C_1 \cap \big\{ ({\bf{v}},{\bf{w}}) \in \Z_{>0}^{6} \ : \ v_1 = 0 \big\}, \\
&C_6 := C_1 \cap \big\{ ({\bf{v}},{\bf{w}}) \in \Z_{>0}^{6} \ : \ w_1 = 0 \big\}.   
\end{align*}

By (\ref{ch. 2, eqn2}), the only relevant values for $c_1 = f H$ are $f = -1,0,1$. The above expression for the generating function gives the following numerical expansions 
\begin{align*}
\sum_{c_{2}} e( \M_{\PP^3}(3,-1,c_{2})) q^{c_{2}} = &3 q^{2} + 42 q^{3} + 333 q^{4} + 1968 q^{5} + 9609 q^{6} + 40881 q^{7} + 156486 q^{8} \\
&+ 550392 q^{9} + 1805283 q^{10} + O(q^{11}), \\
\sum_{c_{2}} e( \M_{\PP^3}(3,0,c_{2})) q^{c_{2}} = &-q^{3} - 9 q^{4} - 60 q^{5} - 309 q^{6} - 1362 q^{7} - 5322 q^{8} - 18957 q^{9} \\
&- 62574 q^{10} + O(q^{11}), \\
\sum_{c_{2}} e( \M_{\PP^2}(3,1,c_{2})) q^{c_{2}} = &3 q^{2} + 42 q^{3} + 333 q^{4} + 1968 q^{5} + 9609 q^{6} + 40881 q^{7} + 156486 q^{8} \\
&+ 550392 q^{9} + 1805283 q^{10} + O(q^{11}).
\end{align*}
This suggests the generating functions $\sum_{c_{2}} e( \M_{\PP^2}(3, \pm c_1,c_{2})) q^{c_{2}}$ are the same. This can be proved by observing that changing $v_{i} \leftrightarrow w_{i}$ and $f \leftrightarrow -f$ swaps terms two $\leftrightarrow$ three and five $\leftrightarrow$ six of the generating function, while leaving terms one and four unchanged. Geometrically, this can be understood as follows. Let $S$ be a nonsingular complete surface, $H$ a polarization, $r > 0$, and $c_{1}$, $c_{2}$ Chern classes. Denote the moduli space of $\mu$-stable locally free sheaves on $S$ of rank $r$ and Chern classes $c_{1}$, $c_{2}$ by $\N_{S}^{H}(r,c_{1},c_{2})$. Then taking the dual $(\cdot)^* = \mathcal{H}{\it{om}}(\cdot, \O_S)$ gives an isomorphism 
\[
\N_{S}^{H}(r,c_{1},c_{2}) \stackrel{\cong}{\longrightarrow} \N_{S}^{H}(r,-c_{1},c_{2}), \ \F \mapsto \F^*.
\]


{\footnotesize{


\end{document}